\documentclass{article}
\usepackage{graphicx}
\usepackage[numbers]{natbib}
\bibpunct{(}{)}{;}{a}{,}{,}
\usepackage[bookmarksnumbered=true, pdfauthor={Wen-Long Jin}]{hyperref}

\newcommand{\commentout}[1]{}

\newcommand{\ba}{\begin{array}}
        \newcommand{\ea}{\end{array}}
\newcommand{\bc}{\begin{center}}
        \newcommand{\ec}{\end{center}}
\newcommand{\bdm}{\begin{displaymath}}
        \newcommand{\edm}{\end{displaymath}}
\newcommand{\bds} {\begin{description}}
        \newcommand{\eds} {\end{description}}
\newcommand{\ben}{\begin{enumerate}}
        \newcommand{\een}{\end{enumerate}}
\newcommand{\beq}{\begin{equation}}
        \newcommand{\eeq}{\end{equation}}
\newcommand{\bfg} {\begin{figure}[h]}
        \newcommand{\efg} {\end{figure}}
\newcommand{\bi} {\begin {itemize}}
        \newcommand{\ei} {\end {itemize}}
\newcommand{\bqn}{\begin{eqnarray}}
        \newcommand{\eqn}{\end{eqnarray}}
\newcommand{\bqs}{\begin{eqnarray*}}
        \newcommand{\eqs}{\end{eqnarray*}}
\newcommand{\bsl} {\begin{slide}[8.8in,6.7in]}
        \newcommand{\esl} {\end{slide}}
\newcommand{\bss} {\begin{slide*}[9.3in,6.7in]}
        \newcommand{\ess} {\end{slide*}}
\newcommand{\btb} {\begin {table}}
        \newcommand{\etb} {\end {table}}
\newcommand{\bvb}{\begin{verbatim}}
        \newcommand{\evb}{\end{verbatim}}

\newcommand{\m}{\mbox}
\newcommand {\der}[2] {{\frac {\m {d} {#1}} {\m{d} {#2}}}}
\newcommand {\pd}[2] {{\frac {\partial {#1}} {\partial {#2}}}}
\newcommand {\pdd}[2] {{\frac {\partial^2 {#1}} {\partial {#2}^2}}}

\newcommand{\mat}[1]{{{\left[ \ba #1 \ea \right]}}}
\newcommand{\cas}[1]{{{\left \{ \ba #1 \ea \right. }}}

\newcommand{\reff}[1] {{{Figure \ref {#1}}}}
\newcommand{\refe}[1] {{{(\ref {#1})}}}
\newcommand{\reft}[1] {{{\textbf{Table} \ref {#1}}}}
\def\l      {{\lambda}}

\def\pmb#1{\setbox0=\hbox{$#1$}%
   \kern-.025em\copy0\kern-\wd0
   \kern.05em\copy0\kern-\wd0
   \kern-.025em\raise.0433em\box0 }

\def\r{{\rho}}
\def\e{{\epsilon}}

\title{A kinematic wave theory of lane-changing traffic flow}
\author{Wen-Long Jin \thanks{Department of Civil and Environmental Engineering, California Institute for Telecommunications and Information Technology, Institute of Transportation Studies, University of California, Irvine, CA 92697-3600. Tel: 949-824-1672. Fax: 949-824-8385. E-mail: wjin@uci.edu. Author for correspondence}}

\begin {document}
\maketitle
\begin{abstract}
Frequent lane-changes in highway merging, diverging, and weaving areas could disrupt traffic flow and, even worse, lead to accidents. 
In this paper, we propose a simple model for studying bottleneck effects of lane-changing traffic and aggregate traffic dynamics of a roadway with lane-changing areas. 
Based on the observation that, when changing its lane, a vehicle affects traffic on both its current and target lanes, we propose to capture such lateral interactions by introducing a new lane-changing intensity variable. 
With a modified fundamental diagram, we are able to study the impacts of lane-changing traffic on overall traffic flow. In addition, the corresponding traffic dynamics can be described with a simple kinematic wave model. For a location-dependent lane-changing intensity variable, we discuss kinematic wave solutions of the Riemann problem of the new model and introduce a supply-demand method for its numerical solutions.
With both theoretical and empirical analysis, we demonstrate that lane-changes could have significant bottleneck effects on overall traffic flow.
In the future, we will be interested in studying lane-changing intensities for different road geometries, locations, on-ramp/off-ramp flows, as well as traffic conditions. The new modeling framework could be helpful for developing ramp metering and other lane management strategies to mitigate the bottleneck effects of lane-changes. 
\end{abstract}


\section{Introduction}
An understanding of the evolution of traffic dynamics is the foundation of transportation network analysis, management, control, and planning. In the celebrated LWR model \citep{lighthill1955lwr,richards1956lwr}, vehicular traffic is viewed as a continuous fluid flow, and traffic dynamics are described by the changes in $(x,t)$-space of three aggregate variables: density $\r$, speed $v$, and flow-rate $q$. The model is can be written as
\bqn
\pd \r t+\pd{\r V(\r)} x&=&0,
\eqn
which is based on three assumptions:
\bi
\item The fundamental law of traffic flow: 
\bqn
q=\r v, \label{fundamentallaw}
\eqn
\item Traffic conservation
\bqn
\pd{\r}{t}+\pd{q}{x}&=&0,\label{trafficconservation}
\eqn
\item Fundamental diagram of speed-density relationship
\bqn
v&=&V(\r). \label{traditionalfd}
\eqn
\ei
The first two assumptions can be derived from the continuum assumption of traffic flow \citep[][Chapters 59, 60]{haberman1977model}. An equilibrium speed-density relationship assumption \refe{traditionalfd} has been shown to be valid in stationary states \citep{delcastillo1995fd_empirical}. Such a relationship can also be derived from various car-following models in steady states \citep{gazis1961follow} and \citep[][Chapter 61]{haberman1977model}. In this sense, the speed-density relationship in \refe{traditionalfd} captures the longitudinal interactions among vehicles, and the traditional LWR model is consistent with car-following behaviors at the aggregate level \citep{newell1961nonlinear}.
In the LWR model, traffic dynamics are described by combinations of shock and rarefaction waves \citep{whitham1974PW}, and the LWR model is thus called the kinematic wave model of traffic flow.

Major roadways, however, usually have multiple lanes, and vehicles can not only change speeds, but also change lanes. That is, vehicles can have both longitudinal and lateral movements.
A lane-changing area is a region, where one or more streams of vehicles systematically change their lanes. Such areas can be near a merging junction and lane-drops,
upstream to a diverging junction, inside a weaving section, or around a cloverleaf interchange
\citep{milam1998cloverleaf,cassidy2005merge}. Since bottlenecks \citep{hall1991capacity} and accidents \citep[e.g.][]{golob2004safety} tend to
occur in these areas, it is important to understand phenomena associated with lane-changing traffic. 
In this study, we assume that different lanes are balanced on average; i.e., vehicle speeds are the same across different lanes at the same location and time. In reality, it has been observed that traffic is nearly balanced inside major weaves,
where speed difference for weaving and non-weaving vehicles is not statistically significant  (about 5 mph)
\citep{roess1974weaving}. Actually, imbalance among different lanes is usually a reason for lane-changing, and lane-changing traffic can have balancing effect; i.e., under certain situations, lane-changes could smooth
out  differences between lanes, and the balancing effect could be beneficial to the whole traffic system in
achieving higher capacity.

At the microscopic level, vehicles' speed adjustment behaviors have been modeled by car-following models \citep{gazis1961follow}; similarly, lane-choice behaviors on why, when, where, and how a vehicle changes its lane are modeled by lane-changing models \citep{gipps1986changing,yang1997dissertation,toledo2003changing,kesting2007lc}. These models can describe detailed lane-changing behaviors, but usually contain a large number of parameters and cannot provide intuitive descriptions of system-level effects of lane-changing traffic. At the macroscopic level, many studies have been carried out to understand various characteristics of lane-changing traffic, including exchange of flows between lanes \citep{michalopoulos1984multilane,holland1997continuum,daganzo2002behavior,Coifman2003inflow}, density oscillation and instability issues \citep{gazis1962multilane,munjal1971multilane1}, and the degree of First-In-First-Out violation among vehicles \citep{jin2006fifo}. In particular, characteristics of weaving sections have been extensively studied since the publication of Highway Capacity Manual in 1950, and many methods have been proposed to analyze levels of service at weaving areas, usually measured by speeds of weaving and non-weaving vehicles  \citep[e.g.][]{leisch1979weaving}. In these methods, parameters  include weave type, geometric configurations, the number of lanes, the length of a weaving area, weaving flow, non-weaving flow, and the number of lane-changes.
In \citep{cassidy1989weaving, cassidy1991weaving, ostrom1993weaving, windover1994weaving}, the distribution of weaving and non-weaving traffic on  rightmost lanes was also studied for different locations.
\citet{eads2000hcm} proposed a framework for evaluating dynamic capacities of a weaving section. 
In \citep{laval2006lc}, a hybrid model of lane-changing traffic was proposed, and it was shown that systematical lane-changes could cause a capacity-drop consistent with observations.
All these studies, however, do not provide a simple approach for analyzing the impacts of lane-changing traffic and the corresponding traffic dynamics at the aggregate level.

In this paper, we attempt to fill this gap and develop a simple kinematic wave model of lane-changing traffic. Based on the observation that, when changing its lane, a vehicle is the leader of two following vehicles on its current and target lanes, we modify the speed-density relationship in \refe{traditionalfd} by adding an effective additional density, equal to total density of all lanes times a lane-changing intensity variable (or simply intensity) $\epsilon(x,t)$. With the new fundamental diagram for both car-following and lane-changing traffic, we then introduce a kinematic wave model based on the LWR model.
When the lane-changing intensity is location-dependent, we analyze the new kinematic wave model as a system of hyperbolic conservation laws. 
We also calibrate a relationship between lane-changing intensity and traffic density with observed vehicle trajectories in a weaving section.
Note that lane-changing intensity $\epsilon(x,t)$ is determined by drivers' lane-changing choices and characteristics in a road section during a time interval. Therefore, $\epsilon(x,t)$ is highly related to the location in a lane-changing section, on-ramp and off-ramp flows, and roads' geometric configurations.
In this paper, we demonstrate that, once $\epsilon(x,t)$ is known, we are able to evaluate the impacts of lane-changing traffic at the aggregate level within the framework of kinematic waves.

The rest of the paper is organized as follows. In Section 2, we introduce a lane-changing intensity variable, a modified fundamental diagram, and a kinematic wave model for lane-changing traffic. In Section 3, we analyze the kinematic wave model as a nonlinear resonant system and propose a corresponding supply-demand method for numerical solutions. In Section 4, we calibrate a relationship between the lane-changing intensity variable and traffic density for a weaving section with NGSIM data.   Finally, some implications and extensions of this study are
discussed.

\section{A kinematic wave model of lane-changing traffic}

\subsection{Model derivation}
\bfg\bc
\includegraphics[width=4in]{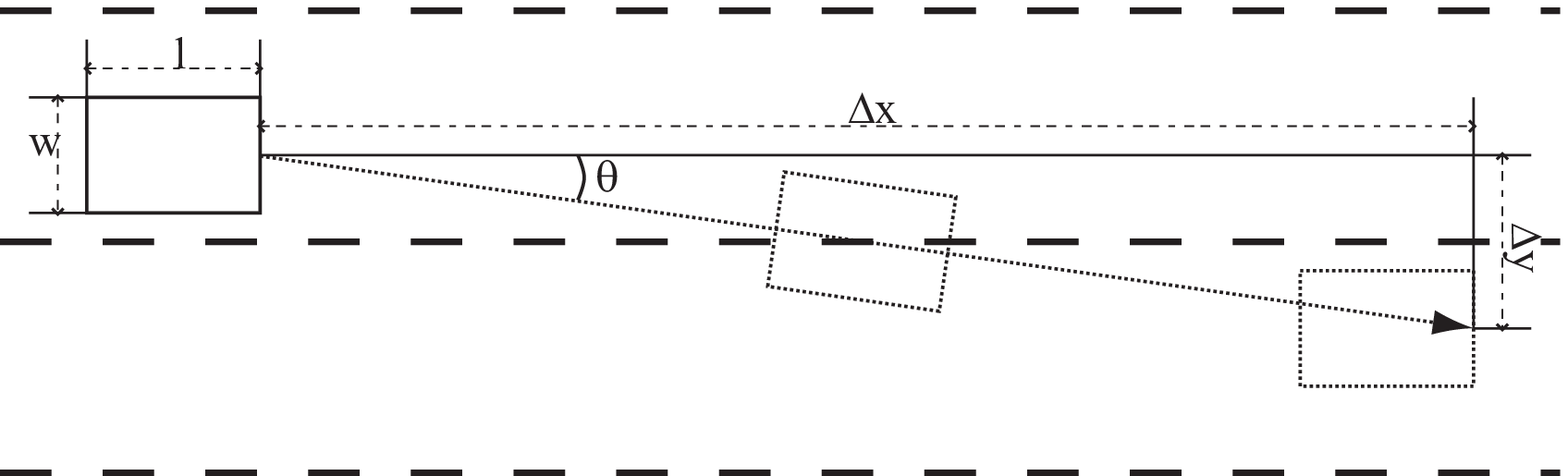} \caption{An illustration of a lane-change} \label{lc20081205illustation}
\ec\efg

Assume the longitudinal axis is $x$, and the lateral axis is $y$. Then the movement of a lane-changing vehicle can be represented by a trajectory in $(x,y,t)$ space. As shown in \reff{lc20081205illustation}, during a time period of $\Delta t$, a lane-changing vehicle can move $\Delta x$ and $\Delta y$ in the $x$ and $y$ directions, respectively. We can see that the lateral displacement threshold of the lane-changing vehicle, or simply lane-changing threshold, $\Delta y$, should be at least the width of the vehicle $w$. 
To accomplish a lane-change, a vehicle first usually signals its intention. In very sparse traffic, it can easily
find a gap big enough to switch to its target lane without waiting too long or interrupting traffic stream on
either its current or target lanes. However, in congested traffic, it usually has to slow down first and wait for
a gap or speed up to squeeze in. In this process, following vehicles on its current lane have to maneuver
accordingly, and following vehicles on the target lane may slow down or switch to other lanes to let it in. That
is, a lane-changing vehicle has longitudinal influence on its current lane and lateral influence on its
subject lane. 
In a lane-changing area, individual vehicles' lane-changing choices may vary significantly \citep{gipps1986changing}, but our experience and many studies suggest that systematic lane-changing vehicles may disrupt traffic flow. 

Since \citep{greenshields1935capacity}, speed- or flow-rate-density relations, $v=V(\rho)$ and $q=\r V(\r)$
respectively, have been used to capture drivers' response to traffic environment at the aggregate level.
Note that, traffic density $\r$ is traditionally defined as the number of vehicles per unit length of road. Thus
traditional fundamental diagram captures only longitudinal interactions between vehicles when they follow each other, but not lateral interactions when they change lanes. Therefore, the traditional LWR model cannot be directly applied to describe traffic dynamics in lane-changing areas.

Since a vehicle uses both its current and target lanes during a lane-change, we propose to double the contribution of a lane-changing vehicle to total density. I.e., lane-changing traffic causes effective additional density. If denoting a lane-changing intensity variable, $\e(x,t)$,  the effective total traffic density is given by
\bqn
\bar \r(x,t)&=&\r(x,t) (1+\e(x,t)). \label{lccoef}
\eqn
Furthermore, for lane-changing traffic, we use the
following modified speed-density relationship:
\bqn
v&=&V(\bar \r)=V(\rho(1+\epsilon)).
\eqn
Then, the fundamental diagram with lane-changing effect is
\bqn
q&=&\r V((1+\epsilon)\r). \label{fd2}
\eqn
Correspondingly, the fundamental diagram without lane-changing effect is
\bqn
q&=&\bar \r V(\bar \r)=(1+\epsilon) \r V((1+\epsilon)\r). \label{fd1}
\eqn
By comparing \refe{fd1} and \refe{fd2}, we will be able to understand the impact of lane-changing traffic.

With the fundamental diagram incorporating lane-changing intensity, we can model lane-changing traffic dynamics in
the framework of kinematic wave theories as follows:
\bqn
\rho_t+\left(\rho V(\rho(1+\epsilon))\right)_x&=&0, \label{kw-weaving}
\eqn
which can be considered as an extension of the LWR model for equilibrium lane-changing traffic.

\subsection{Determination of lane-changing intensity in uniform traffic}\label{uniformtraffic}
Generally $\epsilon(x,t)$ is time- and location-dependent, and $\e(x,t)\geq 0$ due to the bottleneck effect. 
Since $\pd{V}{\epsilon}=V'\r$, which is negative in general, this modified speed-density relationship is consistent with observations \citep{fazio1986weaving} that the total number of lane-shifts in a weaving section, proportional to $\epsilon$, is negatively correlated to weaving and non-weaving speed.
At the microscopic level, lane-changing intensity $\epsilon(x,t)$ is determined by drivers' lane-changing choices and characteristics in a road section during a time interval. At the macroscopic level, the lane-changing intensity variable in \refe{lccoef} are related to locations, road geometric configurations, on-ramp and off-ramp flows, and other traffic conditions. Therefore, in practical applications, $\e(x,t)$ has to be calibrated for different locations of a lane-changing area and traffic conditions.

Here we attempt to determine lane-changing intensity when traffic conditions are uniform in the lane-changing region; i.e., traffic density is the same across the region and for all lanes, and all vehicles travel at the same speed. This is not meant to provide a generic formula for computing $\e$ in \refe{lccoef}. Instead, we intend to show which parameters could affect the lane-changing effect at the aggregate level.

For a lane-changing section of length $L$, we have the
following quantities: $\r_{LC}$ is the density of lane-changing traffic, $\r_{NLC}$ the density of
non-lane-changing traffic, $\r=\r_{LC}+\r_{NLC}$ the total density, $v$ the speed of both lane-changing and non-lane-changing vehicles, $q_{LC}=\r_{LC} v$ the
flow-rate of lane-changing traffic, and $q_{NLC}=\r_{NLC} v$ the flow-rate of non-lane-changing traffic. Then,
the time for all vehicles to traverse the lane-changing region is $T=L/v$. In order to understand how much effect lane-changing traffic can have on total traffic, we have to know the time for finishing a lane-changing maneuver and the number of lane-changes.

We denote the total number of lane-changes in the lane-changing area during a period of $T$ by $N_{LC}$. In literature, the lane-changing frequencies have been studied
for different traffic conditions and road geometric configurations
\citep{Worrall1970lc1,Worrall1970lc,Pahl1972lc,chang1991characteristics,Klar1999multilane}. Here, we are mostly
interested in the relationship between lane-changing traffic flow and the number of lane-changes for relatively
congested traffic.
In \citep{fazio1986weaving}, it has been calibrated that the total number of lane-changes are linear to weaving, i.e., lane-changing, flows for different lane configurations of weaving areas. In \citep{fitzpatrick1996weaving}, it is observed that
a weaving vehicle generally takes 1.33 lane-changes in one-sided weaving operations on one-way frontage roads.
Thus here we simply assume a linear relationship between the number of lane-changes and lane-changing flow; i.e.,
$N_{LC}=\alpha q_{LC} T=\alpha \r_{LC} L$, where the coefficient $\alpha$, the average number of lane-changes of each lane-changing vehicle, could be related to
lane-configurations, number of lanes, traffic conditions, drivers' tendency to lane-changing, and the length of
the section. In addition, $\alpha$ could also depend on location \citep{cassidy1991weaving,windover1994weaving}.

We denote lane-change duration by $t_{LC}$, which starts when a vehicle signals its lane-changing intention and
ends when it finishes a lane-change. If the width of a lane is $D$, and the average lane-changing angle $\theta$, we then
have
\bqn
t_{LC}&=&\frac{D}{v \tan\theta}. \label{lcangle}
\eqn
Note that $\theta$ is related to drivers' behavior, road geometric configurations, and traffic conditions. Usually the time to
complete a lane-changing maneuver is around 2.5 sec on freeways \citep{wangyh1998comparison}, while about 10 sec on
surface streets \citep{sheu2001changing}.

Therefore, by doubling the number of lane-changing vehicles during their lane-changing periods, the effective
number of vehicles at any moment is
\bqn
\bar N&=&\frac{\r  L T+N_{LC} t_{LC}}{T}=\r L+\alpha\r_{LC}L\frac{t_{LC}}{T}.
\eqn
Then we can obtain the lane-changing intensity variable, $\e$, in \refe{lccoef}as
\bqn
\epsilon&=&\frac{N_{LC}t_{LC}}{\r L T}=\frac{N_{LC}t_{LC}}{NT} =\alpha \frac{\r_{LC} t_{LC}}{\r T}. \label{eps}
\eqn
In this sense, the lane-changing intensity can be considered as ratio of the total lane-changing time to the total traveling time during a time interval on a road section.
As expected, lane-changing intensity variable is determined by the number of lane-changes, the lane-changing duration, and traffic density in the lane-changing area. In particular, $\alpha$ and $\r_{LC}$ or $N_{LC}$ can depend on road geometry, locations, and traffic conditions.

\bfg
\bc \includegraphics[height=4cm]{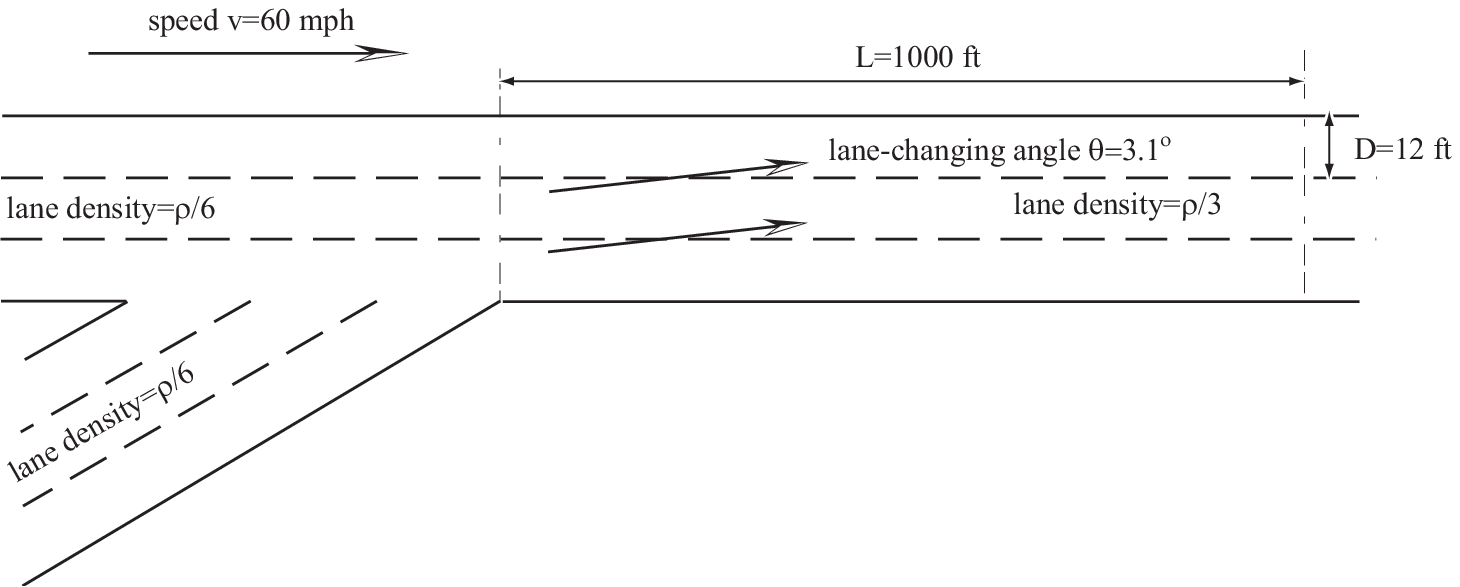}
\caption{A lane-changing area downstream to a merging junction}\label{lcmerge}\ec
\efg

As a simple example, we consider a uniform, three-lane lane-changing region, which is downstream to a merging junction connecting two three-lane roads as shown in \reff{lcmerge} \footnote{Here we simplify the problem by assuming that all lane-changes occur in the region downstream to the merging junction. In reality, some vehicles might change their lanes upstream to the merging junction.}. The length of the lane-changing area is $L=1000$ ft, the width of a lane $D=12$ ft, and both branches have the same density of $\r/ 2$; i.e., traffic density on each lane is $\r/6$ before merging. We assume that all lanes in the lane-changing region are fully balanced; i.e., traffic density on each lane in this region is twice as that before merging, $\rho/3$. Therefore, one third of vehicles from the merging branch, i.e., $\r/6$, have to stay on
the rightmost lane in the lane-changing region with no lane-change, one third will change one lane to the middle lane, and one third will change two lanes to the
leftmost lane. In this case, the total density in the lane-changing region is $\r$, and the density of lane-changing traffic is that from the merging branch; i.e., $\r_{LC}=\r/3$. Then we have the average number of lane-changes of a merging vehicle
\bqs
\alpha=\frac{\r/6\cdot 0+\r/6\cdot 1+\r/6\cdot 2}{\r/3}=1.5.
\eqs
If we assume that vehicle speed $v=$60 mph and the lane-change duration $t_{LC}$=2.5 sec, then the time for a vehicle to traverse the lane-changing region is $T=L/v$=11.4 sec, and we can find from \refe{lcangle} that the lane-changing angle $\theta=3.1$ degrees. Then from \refe{eps}, we can have 
\bqs
N_{LC}t_{LC}&=&\frac{\r}6 L t_{LC}+\frac{\r}6 L 2t_{LC}=\frac{\r}2 Lt_{LC},
\eqs 
and the lane-changing intensity variable in \refe{eps} is
\bqs
\epsilon=\frac 12 \frac{2.5}{11.4}=0.11.
\eqs
Here we assume that lane-changes are evenly distributed over the whole road section. In reality, however, the distribution may be uneven, and we may have a location-dependent lane-changing intensity variable $\epsilon(x)$.

\subsection{Bottleneck effects of lane-changing traffic}
We consider the following triangular fundamental diagram
\citep{munjal1971multilane,haberman1977model,newell1993sim},
\bqs
V(\r)&=&\cas{{ll}v_f, & 0\leq \r\leq \r_c;\\ \frac{\r_c}{\r_j-\r_c} v_f \frac{\r_j-\r}{\r} , &\r_c< \r\leq \r_j,
}\\
Q(\r)&=&\cas{{ll}v_f \r, & 0\leq \r\leq \r_c;\\ \frac{\r_c}{\r_j-\r_c} v_f (\r_j-\r) , &\r_c< \r\leq \r_j,  }
\eqs
where $v_f$ is the free flow speed, $\r_j$ the jam density, and $\r_c$ the critical density where flow-rate,
$q=\r v$, attains its maximum, i.e. the capacity. The values of these parameters are: $v_f=65$ mph, $\r_j=240$
vpmpl, $\r_c= \r_j/6=40$ vpmpl, and the capacity is $Q_{max}=2600$ vphpl \citep{delcastillo1995fd_empirical}.
After introducing lane-changing effect, we then have the following fundamental diagram,
\bqn
V(\e,\r)&=&\cas{{ll}v_f, & 0\leq \r\leq \r_c/(1+\epsilon);\\ \frac{\r_c}{\r_j-\r_c} v_f
\frac{\r_j-\r(1+\epsilon)}{\r(1+\epsilon)} , &\r_c/(1+\epsilon)< \r\leq \r_j/(1+\epsilon),  }\\
Q(\e,\r)&=&\cas{{ll}v_f \r, & 0\leq \r\leq \r_c/(1+\epsilon);\\ \frac{\r_c}{\r_j-\r_c} v_f
\frac{\r_j-\r(1+\epsilon)}{1+\epsilon} , &\r_c/(1+\epsilon)< \r\leq \r_j/(1+\epsilon). }\label{fd_tri}
\eqn

\bfg
\bc \includegraphics[height=8cm]{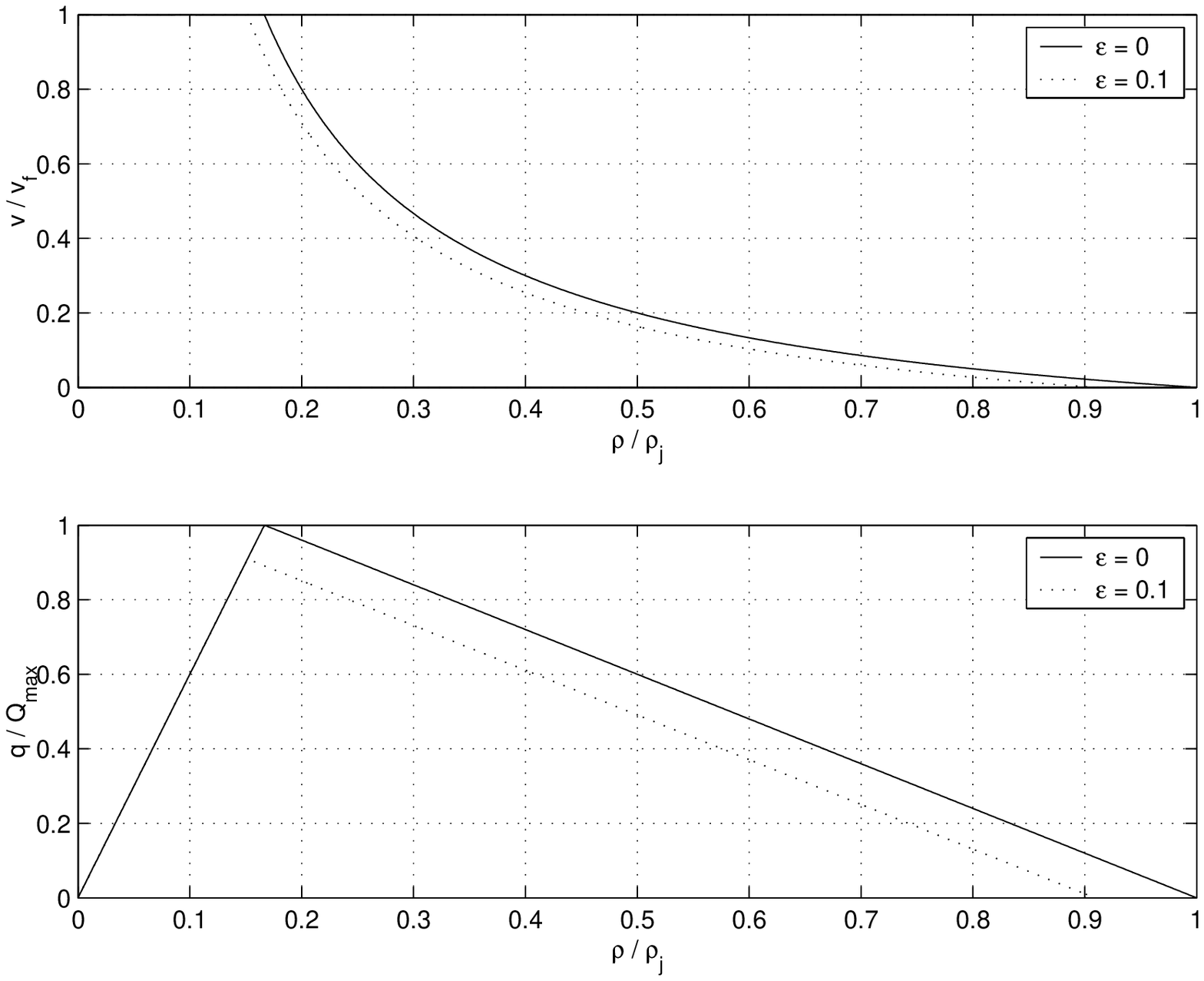}
\caption{The triangular fundamental diagram with constant lane-changing effect}\label{fdconstant}\ec
\efg

For one example, we assume constant lane-changing effect, $\epsilon=0.1$. The relationships between the
standardized density, speed, and flow-rate are shown in \reff{fdconstant}. From the figure, we can see that: (i)
lateral interactions can be omitted when traffic is relative sparse; (ii) there exists capacity reduction of
$1-1/1.1$=9.1\%, which is consistent in magnitude with observations in
\citep{cassidy1999bottlenecks}\footnote{Note that reasons of bottleneck were not discussed in this study, and we
suspect that systematic lane-changes caused by merging traffic from on-ramps could have significant
contributions.}; and (iii) the observed jam density is lower than maximum jam density, which could be one reason
of observing different jam densities \citep{delcastillo1995fd_empirical}.
That is, for the same road section, constant lane-changes can yield a lower throughput. This clearly demonstrates the bottleneck effect of lane-changing traffic.

For another example, we consider the following density-dependent lane-changing effect, which has a jump at critical density,
\bqn
\epsilon(\r)&=&\cas{{ll}0, &\r<\r_c;\\
\frac{2-2\r/\r_j}{15+2\r/\r_j},&\r\geq \r_c.}
\eqn
As shown in \reff{fd20050119rlambda}, the lane-changing effect is approximately linear for congested traffic, and
the resulted fundamental diagram has a shape of reverse-$\lambda$ \citep{koshi1983fd}.
Since lane-changing intensity is related to the location, road geometry, on-ramp and off-ramp flows, and other exogenous conditions, it is possible to have a discontinuous $\r-\epsilon$ relationship. In \citep{cassidy2005merge}, for example, it was observed that lane-changing intensity is significantly higher when a 16-vehicle queue appears on the shoulder lane.
Therefore, this discontinuous $\r-\epsilon$ relationship could also be a cause of capacity-drop: when traffic density increases over a critical value, the number of lane-changes significantly increases, and the overall capacity drops as a result.

\bfg
\bc \includegraphics[height=8cm]{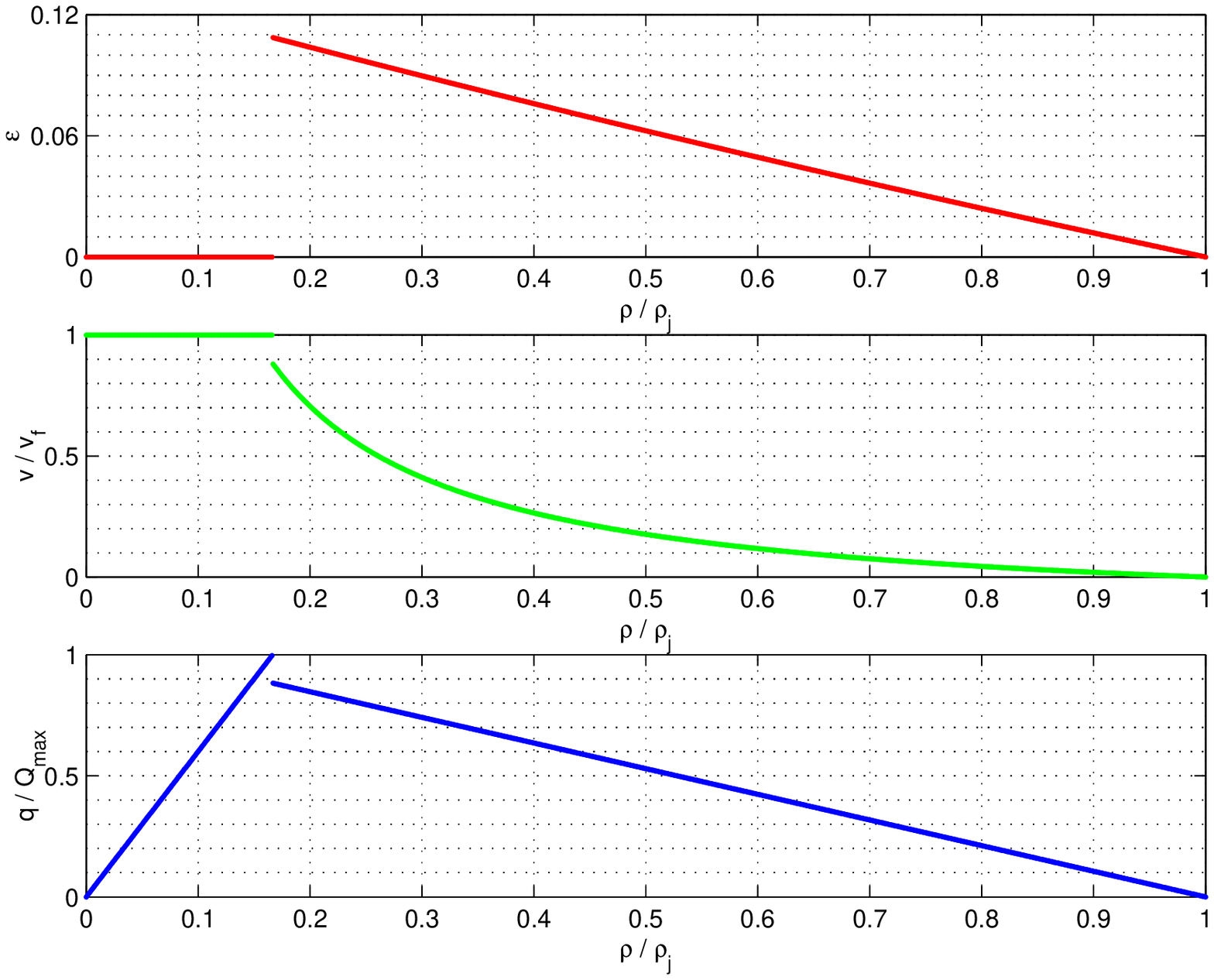}
\caption{A fundamental diagram of reverse-$\lambda$ shape with density-dependent lane-changing
effect}\label{fd20050119rlambda}\ec
\efg

\section{Traffic dynamics for location-dependent lane-changing intensities}
In this section, we consider a special case, where the lane-changing intensity variable, $\epsilon(x,t)$, is independent of
traffic conditions and only a function of location denoted by $\epsilon(x)$. Although simplified, this assumption does capture the most pronounced feature that lane-changing intensity is location-dependent. For this scenario, we can analyze \refe{kw-weaving} and understand some fundamental traffic dynamics in lane-changing areas. In the literature, there have been many methods for analyzing such an inhomogeneous LWR model with discontinuous flux functions (refer to \citep{jin2009sd} for a survey). In this study, we follow an approach proposed in \citep{jin2003inhlwr}.

\subsection{A system of hyperbolic conservation laws}
For a location dependent lane-changing intensity $\epsilon(x)$, we obtain a simple conservation law of lane-changing coefficient,
\bqn
\epsilon_t&=&0,
\eqn
which is equivalent to saying that $\epsilon$ is time-independent.

Without considering inhomogeneities of links \citep{jin2003inhlwr} or merging and diverging effect, we can have the following system of conservation laws,
\bqn
U_t+F(U)_x&=&0, \label{hcl-weaving}
\eqn
where $U=(\e,\r)$, $F(U)=(0, Q(\e,\r))$, $x\in R$, and $t\geq 0$. Note that the domains of $\e$ and $\r$  are
$\e\in [0,\infty)$ and $\r\in[0,\frac{\r_j}{1+\e}]$ respectively.

Here we use a differentiable, concave fundamental diagram in our analysis. For example, we can use the following
maximum sensitivity model  \citep{delcastillo1995fd_empirical}:
\bqn
V(\r)&=&v_f\left\{ 1-\exp\left[1-\exp\left( \frac{|c_j|}{v_f}(\frac{\r_j}{\r}-1)\right)\right]\right\},
\eqn
where $c_j$ is the shock wave speed for jammed traffic.
In this fundamental diagram, $V'<0$ and $\frac{\m{d}^2}{\m{d}\r^2} Q(\r)<0$.

For the kinematic wave model of lane-changing traffic, \refe{hcl-weaving}, the second term of flux, $F(U)_x$, can
be linearized as $\partial F(U) U_x$ with
\bqs
\partial F&=&\mat{{cc} 0&0\\\r^2 V' & V+\r (1+\e) V'},
\eqs
whose two eigenvalues, or wave speeds, are $\l_0(U)=0$, and $\l_1(U)=V+\r \e V'=Q'((1+\e)\r)$.\footnote{By
assuming $\bar \r=\r/(1+\e)$, we have that $Q=\frac{1}{1+\e}\bar \r V(\bar \r)=\frac{1}{1+\e} Q(\bar \r)$ and
$\l_1(U)=\pd{Q}{\r}=\frac{1}{1+\e}\der{Q}{\bar \r} \pd{\bar \r}{\r}$. Therefore, $\l_1(U)=\der{Q}{\bar
\r}=Q'((1+\e)\r)$.} The corresponding eigenvectors are
\bqs
\vec R_0&=&\mat{{c}V+\r (1+\e) V'\\-\r^2 V'}, \quad \vec R_1=\mat{{c} 0 \\1}.
\eqs
 Since $V\geq 0$ and $V'< 0$, it is possible that $\l_1=0=\l_0$. Therefore, \refe{hcl-weaving} is a non-strictly
hyperbolic system.

In addition, if defining the critical traffic state $U_*=(\e_*,\r_*)$ by $\l_1(U_*)=0$, we then have the
following results.
First,
\bqn
\pd{}{\r}\l_1(U_*)&=&\pdd{Q}{\r}|_{U_*}=(1+\e)Q''((1+\e)\r)|_{U_*}<0, \label{decreasingev}
\eqn
since $Q$ is strictly concave in $(1+\e)\r$.
Second,
\bqs
\pd{}{\e} Q(U_*)&=&\r^2 V'|_{U_*}\leq 0.
\eqs
From \refe{decreasingev}, the wave speed $\l_1$ is decreasing with respect to traffic density. In addition,
$\l_1(\e,\r=0)=v_f>0$ and $\l_1(\e,\r)=\e\r V'<0$ when $(1+\e)\r=\r_j$. Therefore, for any weaving factor $\e_*$,
we can find a unique critical state $U_*=(\e_*,\r_*)$.  Further, since $\partial\l_1/\partial \rho<0$, $Q(U_*)$
is indeed the maximum flow-rate for $\e_*$ and equals $\frac{1}{1+\e_*}\max\{Q(\e,\r)|\e=0\}$, where
$\max\{Q(\e,\r)|\e=0\}$ is the capacity without lane-changes, and $\frac 1 {1+\e_*}$ reflects capacity-drop caused by
lane-changes. Here we call $\r_*=\Gamma(\e_*)$ as the transitional curve, since traffic states on its left side
are under-critical (UC) while those on its right side are over-critical (OC). From these properties,
\refe{hcl-weaving} is a nonlinear resonant system and, in the neighborhood of a critical state, the Riemann
problem can be uniquely solved \citep{isaacson1992resonance}.

\subsection{Solutions of the Riemann problem}
We consider the  Riemann problem of \refe{hcl-weaving} with the following initial conditions,
\bqs
U(x,t=0)&=&\cas{{ll} U_L, &x<0, \\U_R, &x>0,}
\eqs
where $U_L=(\e_L,\r_L)$ and $U_R=(\e_R,\r_R)$. Since, in the Godunov method \citep{godunov1959}, general initial
conditions can be approximated by piece-wisely linear functions, we can solve the Riemann problem at each
boundary to find the flux and then update traffic conditions for the next time step.

\bfg
\bc \includegraphics[height=8cm]{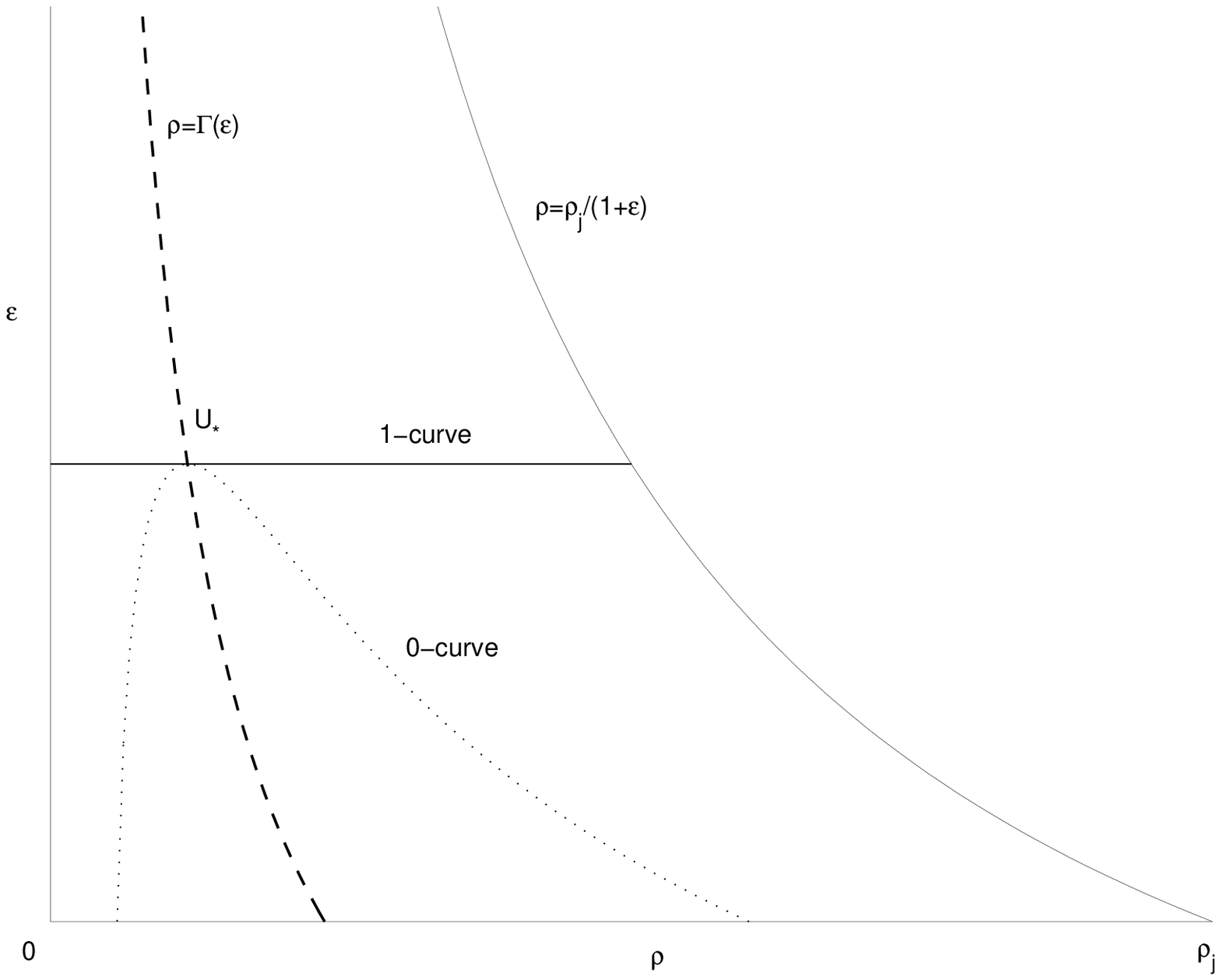}
\caption{A 0-curve, 1-curve, and the transitional curve in $(\r,\e)$-plane}\label{phasecurve}\ec
\efg

According to \citep{isaacson1992resonance}, the Riemann problem is solved by a combination of two basic types of
waves, which are associated with the two eigenvalues of $\partial F$:  we have 0- or
standing waves with wave speed $\l_0=0$, and 1-waves with wave speed $\l_1$. Note that 0-waves carry contact discontinuities, and 1-waves include traditional shock and rarefaction
waves. In the $\r-\e$ phase plane, $(U_L,U_R)$ yielding 0-waves forms a 0-curve, and similarly we can define a
1-curve. Further, 0-curves are the integral curves of eigenvector $\vec R_0$, and 1-curves are those of $\vec
R_1$. They are $\r V(\r(1+\e))=$const and  $\e=$const, respectively. A 0-curve, a 1-curve, and the transitional
curve are shown in \reff{phasecurve}. For a Riemann problem, its solutions must satisfy both Lax's
\citep{lax1972shock} and Isaacson and Temple's \citep{isaacson1992resonance} entropy conditions, such that
resulted waves increase their wave speeds from left to right, and 0-curves do not cross the transitional
curve. Since both 0- and 1-waves are linear in the sense that $U(x,t)=U(x/t)$, then from the Riemann solutions, we can
obtain the boundary flux through $x=0$ as $q=Q(U(x=0,t>0))$.

\bfg
\bc \includegraphics[height=8cm]{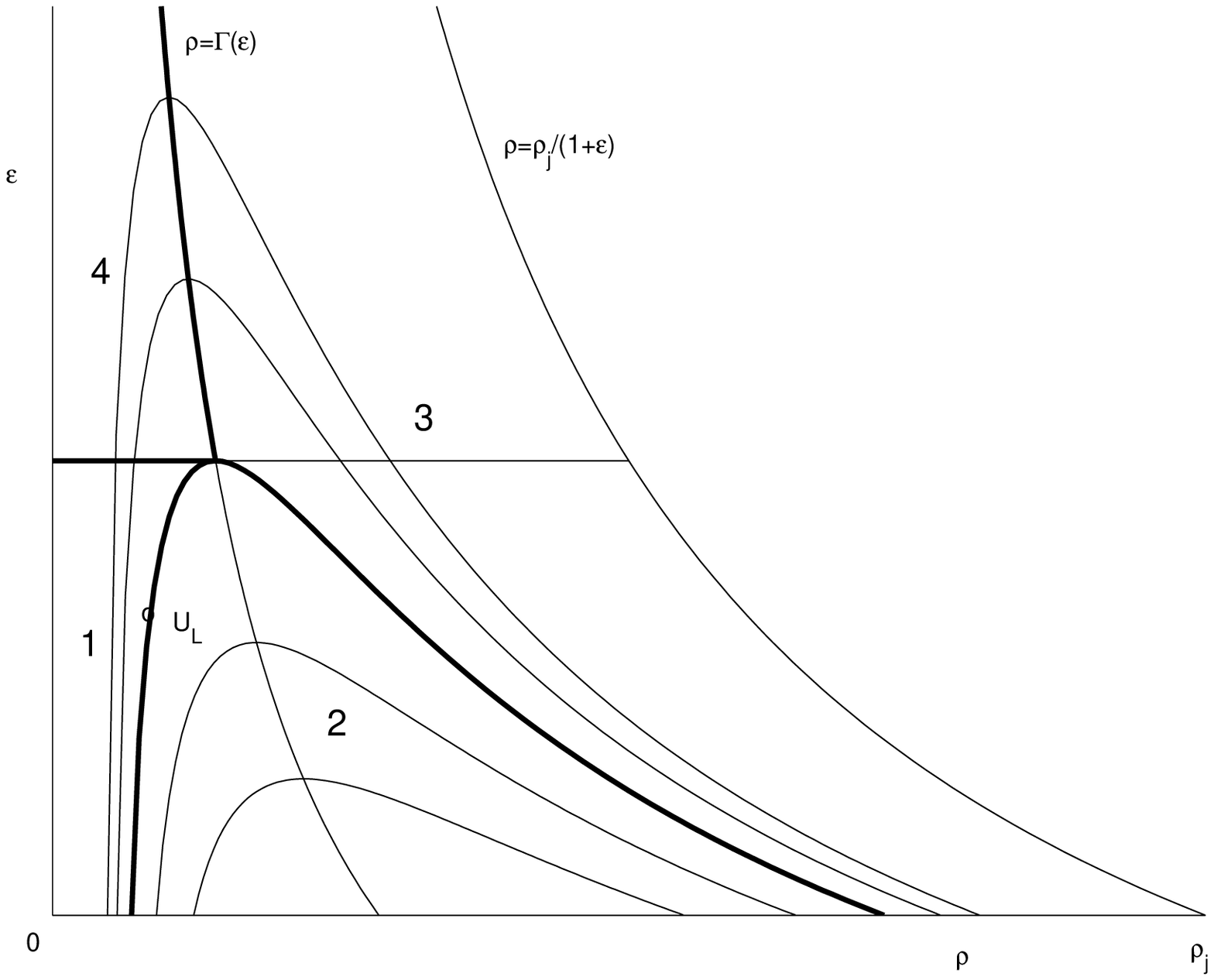}
\caption{Four types of Riemann solutions when $U_L$ is UC}\label{riemannsolutionsL}\ec
\efg

As shown in \reff{riemannsolutionsL}, when $U_L$ is UC, i.e. to the left of the transitional curve in the phase
plane, we have $\r_L\leq \Gamma(\e_L)$ and can divide the $\r-\e$ phase plane into four regions by the thick curves and obtain the four
types of wave solutions for different values of $U_R$, which satisfy the two aforementioned entropy conditions. Here kinematic wave solutions of Type 1 are carefully explained, and solutions of Types 2 to 4 are given in Appendix.
\bi
\item[Type 1] When $U_R$ is in Region 1, $\r_R\leq \Gamma(\e_R)$, $Q(U_R)\leq Q(U_L)$, and $Q(\e_R, \Gamma(\e_R))\geq Q(U_L)$. That is, $U_R$ is UC, $Q(U_R)$ is not greater than $Q(U_L)$, but the capacity at $\e_R$ is not smaller than $Q(U_L)$. The Riemann problem is solved by a combination of a standing wave and a
forward traveling rarefaction wave, with an intermediate state $U_1$ as shown by the four figures of
\reff{lctype1}. The top left figure in the $\r-\e$ phase plane shows that $(U_L,U_1)$ is a standing wave, and the
two states share the same flow-rate $Q(U_L)$. The bottom left figure shows two fundamental diagrams for when
$\e=\e_L$ and $\e=\e_R$ respectively and the Riemann solution of a standing wave and a rarefaction wave. The top
right figure shows characteristic waves in the $x-t$ plane, and we can see the discontinuity caused by the
standing wave at $x=0$ and the rarefaction wave on its right side. The bottom right figure shows the profile of
density at a time instant $t_0$. In this case, since the flow-rate remains the same across $x=0$, then the
boundary flux $q=Q(U_L)$.

\bfg
\bc \includegraphics[height=8cm]{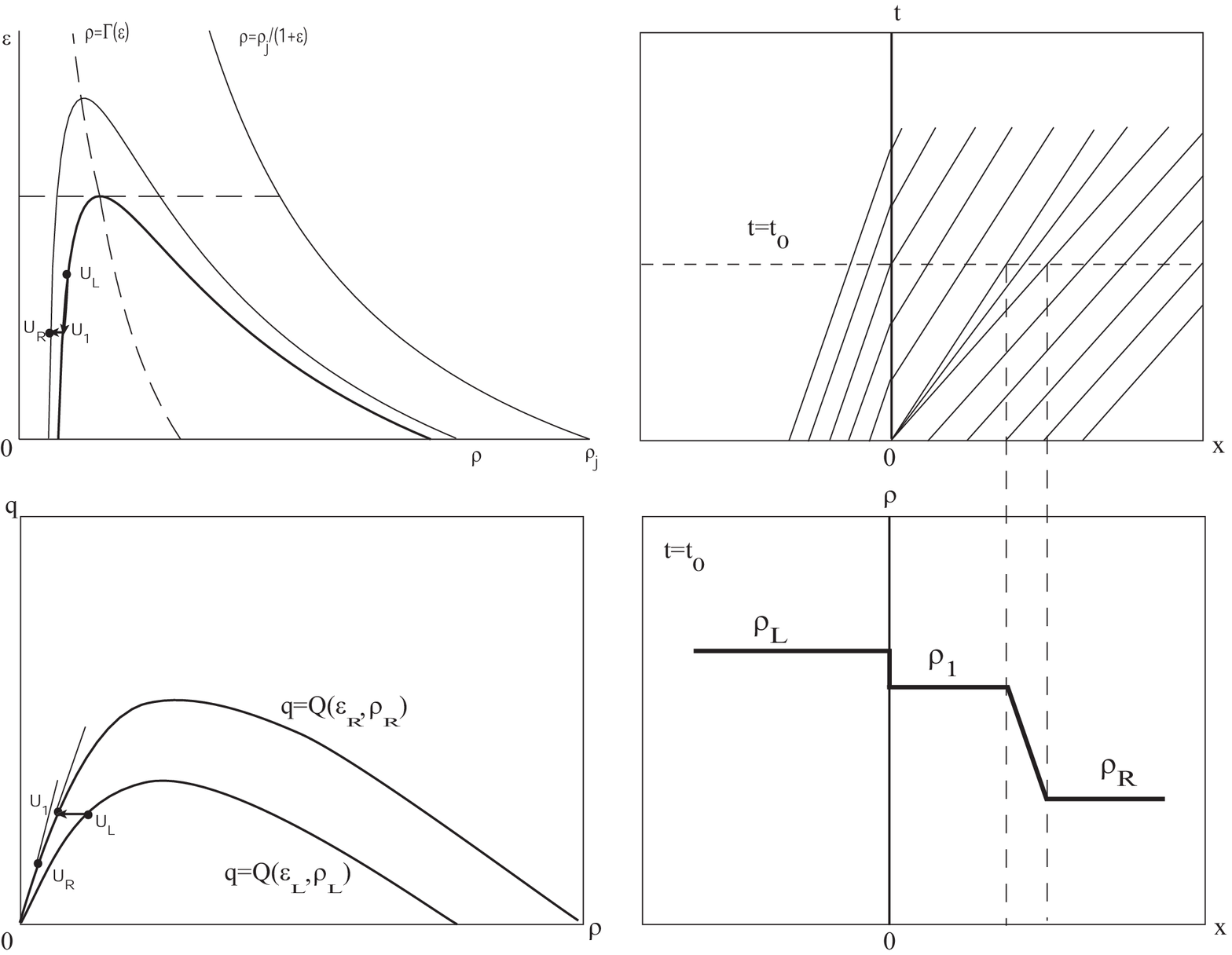}
\caption{An example for wave solutions to Riemann problem of Type 1}\label{lctype1}\ec
\efg

\ei

\bfg
\bc \includegraphics[height=8cm]{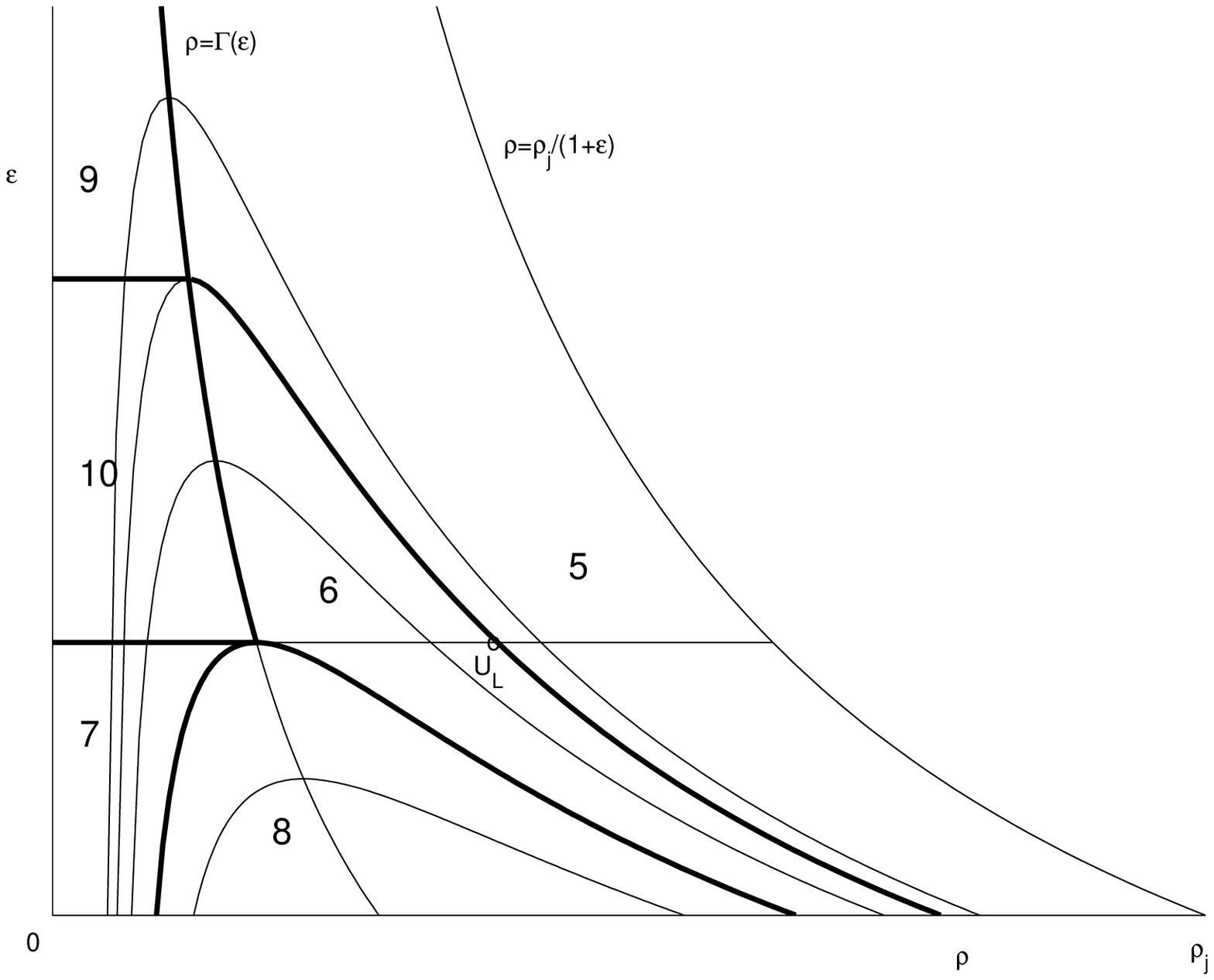}
\caption{Six types of Riemann solutions when $U_L$ is OC}\label{riemannsolutionsR}\ec
\efg
As shown in \reff{riemannsolutionsR}, when $U_L$ is OC, i.e. to the right of the transitional curve in the phase
plane, we have $\r_L> \Gamma(\e_L)$ and can divide the $\r-\e$ phase plane into six regions by the thick curves and obtain the six
types of wave solutions for different values of $U_R$, which satisfy the two aforementioned entropy conditions. The solutions of Types 5 to 10 are also given in Appendix.

If we introduce new definitions of local traffic supply and demand \citep{daganzo1995ctm, lebacque1996godunov,
nelson2004model} as follows,
\bqn
\ba{lcl}
S(\e,\r)&=&\max\{Q(\e,\bar\r): \r\leq \bar\r\leq \r_j\},\\
D(\e,\r)&=&\max\{Q(\e,\bar\r): 0\leq\bar\r\leq\r\}.
\ea\label{def:sd}
\eqn
This is equivalent to saying that
\bqs
S(\e,\r)&=&\cas{{ll} Q(\e,\Gamma(\e)), &\r\in[0, \Gamma(\e)],\\Q(\e,\r), &\r\in(\Gamma(\e), \r_j];}
\eqs
and
\bqs
D(\e,\r)&=&\cas{{ll} Q(\e,\r), &\r\in[0, \Gamma(\e)],\\Q(\e,\Gamma(\e)), &\r\in(\Gamma(\e), \r_j],}
\eqs
where $Q(\e,\Gamma(\e))$ is the capacity at $\e$.
These definitions lead to that both local supply and demand are not smaller than flow-rate; i.e., $S(\e,\r)\geq Q(\e,\r)$, and $D(\e,\r)\geq Q(\e,\r)$.

\btb
\bc
\begin{tabular}{|@{}c@{}||c@{}|c@{}||c@{}|c@{}||c@{}|}\hline
Type & $U_L$ &$D(U_L)$&$U_R$&$S(U_R)$&$q$\\\hline
1& UC & $Q(U_L)$ & UC, $Q(U_R)\leq Q(U_L)\leq Q(\e_R,\Gamma(\e_R))$ &$Q(\e_R,\Gamma(\e_R))$&$Q(U_L)$\\\cline{1-1}\cline{4-6}
2&&&$Q(U_R)>Q(U_L)$&$\geq Q(U_R)$&$Q(U_L)$\\\cline{1-1}\cline{4-6}
3&&&OC, $Q(U_R)\leq Q(U_L)$&$Q(U_R)$&$Q(U_R)$\\\cline{1-1}\cline{4-6}
4&&&UC, $Q(\e_R,\Gamma(\e_R))<Q(U_L)$&$Q(\e_R,\Gamma(\e_R))$&$Q(\e_R,\Gamma(\e_R))$\\\hline
5&OC&$Q(\e_L, \Gamma(\e_L))$ &OC, $Q(U_R)\leq Q(U_L)$&$Q(U_R)$&$Q(U_R)$\\\cline{1-1}\cline{4-6}
6&&$\geq Q(U_L)$&OC, $Q(U_L)< Q(U_R)\leq Q(\e_L, \Gamma(\e_L))$&$Q(U_R)$&$Q(U_R)$\\\cline{1-1}\cline{4-6}
7&&&UC, $Q(U_R)\leq Q(\e_L, \Gamma(\e_L))\leq Q(\e_R, \Gamma(\e_R))$&$Q(\e_R, \Gamma(\e_R))$&$Q(\e_L, \Gamma(\e_L))$\\\cline{1-1}\cline{4-6}
8&&&$Q(U_R)>Q(\e_L, \Gamma(\e_L))$&$\geq Q(U_R)$&$Q(\e_L, \Gamma(\e_L))$\\\cline{1-1}\cline{4-6}
9&&&UC, $Q(\e_R, \Gamma(\e_R))<Q(U_L)$&$Q(\e_R, \Gamma(\e_R))$&$Q(\e_R, \Gamma(\e_R))$\\\cline{1-1}\cline{4-6}
10&&&UC, $Q(U_L)\leq Q(\e_R, \Gamma(\e_R))<Q(\e_L, \Gamma(\e_L))$&$Q(\e_R, \Gamma(\e_R))$&$Q(\e_R, \Gamma(\e_R))$\\\hline
\end{tabular}
\caption{Boundary flux, downstream supply, and upstream demand for ten types of Riemann problems}\label{tcomp}\ec
\etb

In \reft{tcomp}, we summarize the solutions of boundary flux for ten types of Riemann problems discussed in the preceding subsection and compute the corresponding upstream demand and downstream supply functions. From the table, we can conclude that the boundary flux can be computed by the minimum of upstream demand and downstream supply; i.e.,
\bqn
q&=&\min\{S(U_R),D(U_L)\}. \label{s-d}
\eqn
Although \refe{s-d} is proved for time-independent and differential fundamental diagram, we expect that  the supply-demand definitions in \refe{def:sd} can also be applied to the triangular fundamental
diagram and to the scenarios when $\e$ is time-varying.
Therefore, this method can be considered as a simplification and extension to the Riemann solver in the preceding subsection and is more
efficient for numerical simulations.

\section{Calibration of lane-changing intensity in a weaving section}
In this section, we consider lane-changing traffic on a multi-lane freeway section on interstate 80 in Emeryville (San Francisco), California, as shown in \reff{i80studyarea}. The freeway section has six lanes, where lane 1 is a car-pool lane, and an on-ramp from Powell Street and an off-ramp to Ashby Ave. The width of lanes 1 to 6 is 12 ft. Traffic direction is from south to north. The road section is covered by seven cameras, numbered 1 to 7 from south to north, and vehicle trajectories were transcribed from video data by FHWA's NGSIM project \citep{fhwa2006ngsim}. In total there are four data sets \citep{ngsim2004data235,ngsim2005data400,ngsim2005data500,ngsim2005data515}: Data set 1 contains trajectories of all vehicles every fifteenth second on December 3, 2003 between 2:35pm and 3:05pm; Data sets 2 to 4 contain the locations of each vehicle every tenth second on April 13, 2005 between 4pm and 4:15pm, between 5pm and 5:15pm, and between 5:15pm and 5:30pm, respectively.  Note that the lengths of the study location in the first and the other three data sets are different: The Ashby Ave off-ramp is included in the first data set, but not in the other three. 

\bfg\bc
\includegraphics[width=6in]{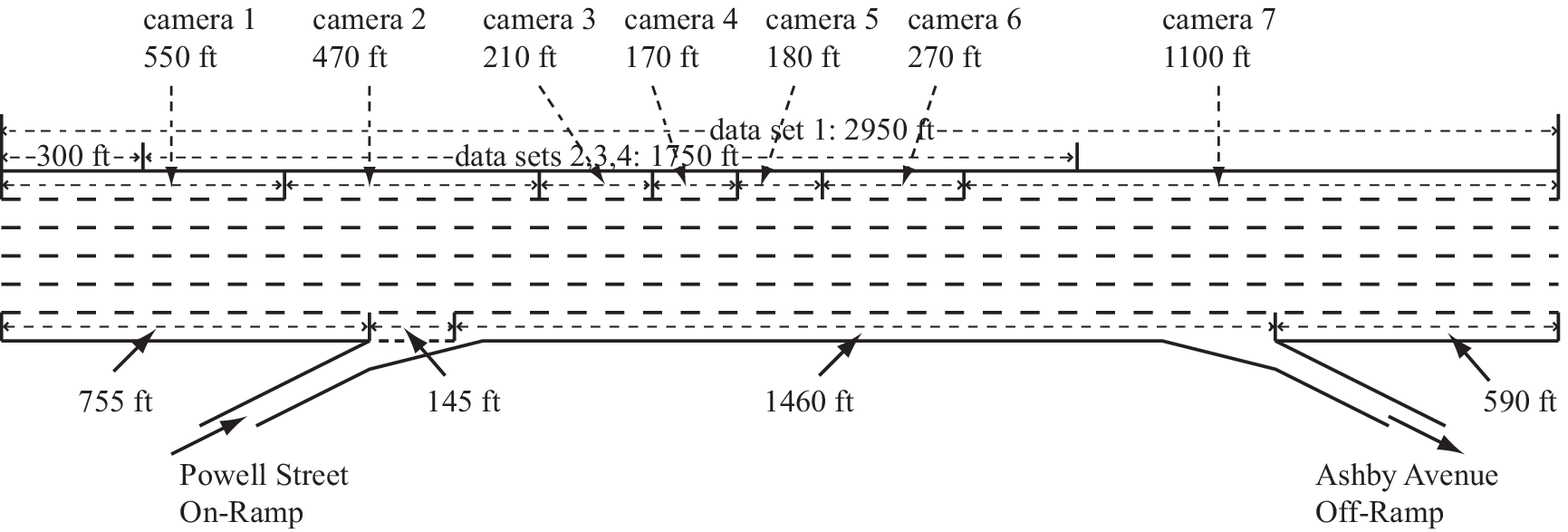}\caption{I-80 study area}\label{i80studyarea}
\ec\efg

\bfg\bc
\includegraphics[width=4in]{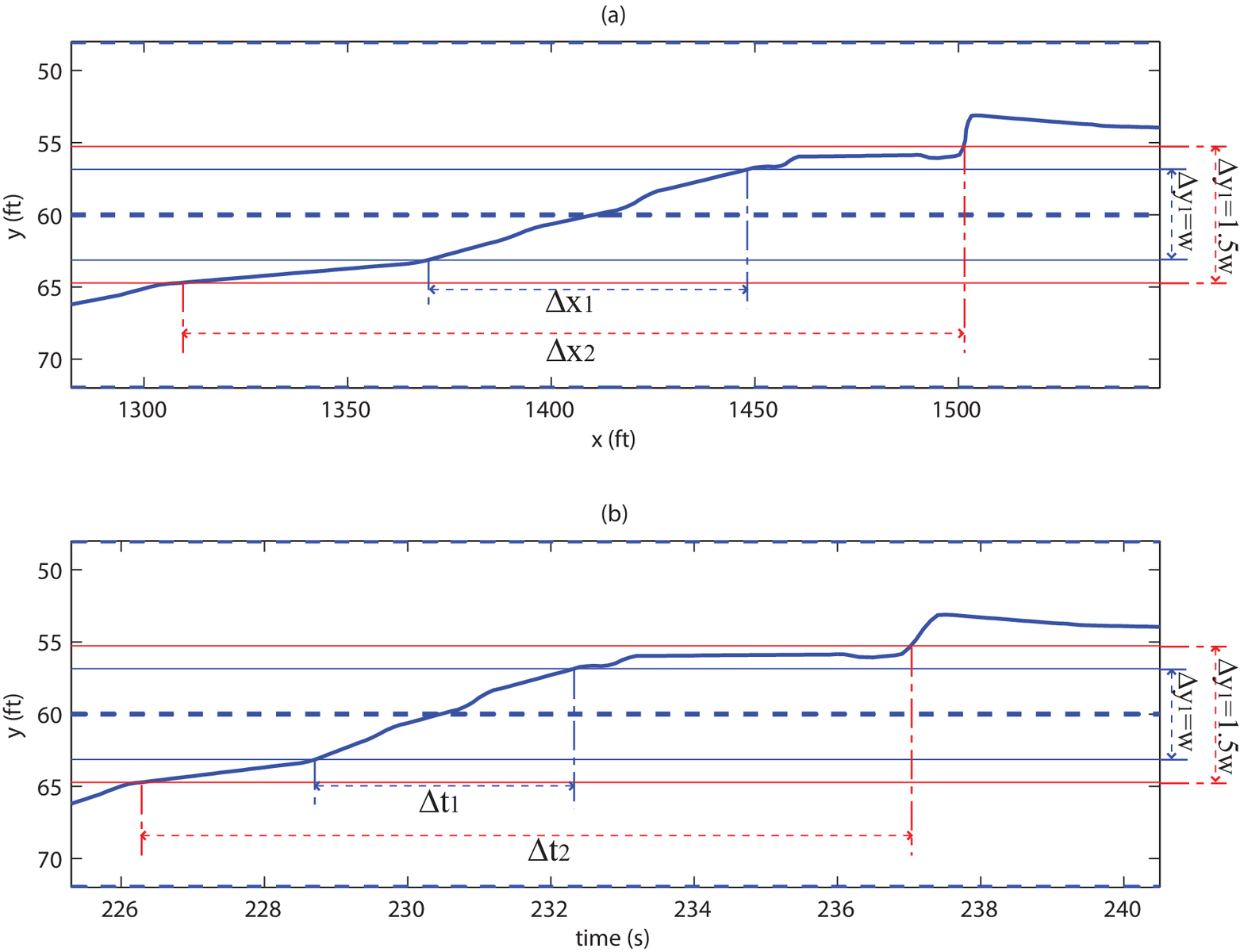}\caption{Trajectories of a lane-changing vehicle: (a) in $x-y$ space; (b) in $t-y$ space}\label{lc20081205trajectory}
\ec\efg

In \reff{lc20081205trajectory}, we demonstrate the trajectories in both $x$-$y$ and $t$-$y$ spaces of vehicle 576 in data set 2, when it switches from lane 6 to lane 5. Here $y=60$ is the lane separation line, and the width of the vehicle is $w=$6.3 ft. In both figures, we demonstrate $\Delta x$ and $\Delta t$ for two different $\Delta y$, $w$ or $1.5w$. From both figures, we can see that the lane-changing angles is decreased by more than half and the lane-changing time is more than doubled when we increase $\Delta y$ from $w$ to $1.5w$. Note that $1.5w$ is still smaller than the width of a lane, 12 ft. From the figures, we can see that, after entering the target lane, the vehicle first stay near the lane separation line and then adjust to the center. Such behavior can be confirmed with video data by camera 6.

Following the example in Section 2.2, we can theoretically derive the lane-changing intensity variable in \refe{eps} and corresponding fundamental diagrams in \refe{fd1} and \refe{fd2} for a lane-changing area with a width of $L=x_b-x_a$ and time period with a duration of $T=t_b-t_a$, shown in \reff{i80studyarea}.
Assume that all lane-changes are induced by on-ramp vehicles and finished in a region of $L$ downstream to the merge and that traffic flow is uniform with speed $v$ and flow-rate $q$, respectively. Following the example in Section \ref{uniformtraffic}, we can estimate the number of lane-changes in $L$ during $T=L/v$ as $2.5 q_{on}T$, since one sixth of the on-ramp vehicles will make zero up to five lane-changes. Thus the lane-changing intensity variable can be estimated as
\bqn
\epsilon&=&\frac{N_{LC}t_{LC}}{\r L T}=\frac{2.5 q_{on} T t_{LC}}{\r L T}=2.5 q_{on} \frac{ t_{LC} }{\r L}. \label{approxe}
\eqn
Since $q_{on}$ is usually controlled by ramp meters and $t_{LC}$ is almost constant, $\epsilon$ is generally not a constant in $\r$ or $v$: it is decreasing in $\r$ and increasing in $v$, when traffic is congested. In particular, if $q_{on}=800$ vph, $\r=200$ vpm, $L=900$ ft, and $t_{LC}=5$ s, we have $\epsilon=8.15\%$ for $v=60$ mph.
Note that \refe{approxe} is just an estimation for the weaving region in \reff{i80studyarea}, since (1) lane-changes can also be caused by vehicles leaving to the off-ramp or vehicles from or to other adjacent on- and off-ramps, (2) lane-changes induced by the on-ramp vehicles can occur upstream to the merge, and (3) fewer than one-sixth of the on-ramp vehicles can make five lane-changes, since the first lane is a car-pool lane.

In the following, we consider a road section of 900 ft downstream to the merge in \reff{i80studyarea} with $x_a=950$ ft and $x_b=$1850 ft. The length is $L=$900 ft or 274 m. 
The time interval is $[t_i, t_i+T]$, where the smallest $t_i$ is the time when the first vehicles pass $x_b$, and the greatest $t_i+T$ is the time when the last vehicles pass $x_a$. For four data sets, $T$=39, 63, 86, 100 s, respectively.
Here we consider different lane-changing thresholds $\Delta y$. Refer to \citep{jin2010lc} for detailed methods for computing $\theta$, $\epsilon$, $\r$, $q$, and $v$.

With a lane-changing threshold of $\Delta y=w$, we have the following results.
From \reff{lc20081207characteristics}, we find a linear relationship between density and lane-changing angle
\bqs
\theta&=&-0.5016+0.0121 \r,
\eqs
with R-square 0.9441.
Here we can also calibrate a relationship similar to \refe{approxe} as 
\bqs
\epsilon&=&-0.0247+\frac{24.2712}{\r},
\eqs
with R-square 0.8077. Hereafter we only calibrate the exponential function, since it usually yields better R-square.
From \reff{lc20081207fd}, we can see that $q=\r V(\r)$ has a capacity of 13408 mph when $\r=233$ vpm, and $q=\r V((1+\epsilon) \r)$ has a capacity of 12369 mph when $\r=215$ vpm. Lane-changes can cause $7.75\%$ reduction in capacity. 

\bfg\bc
\includegraphics[height=4in]{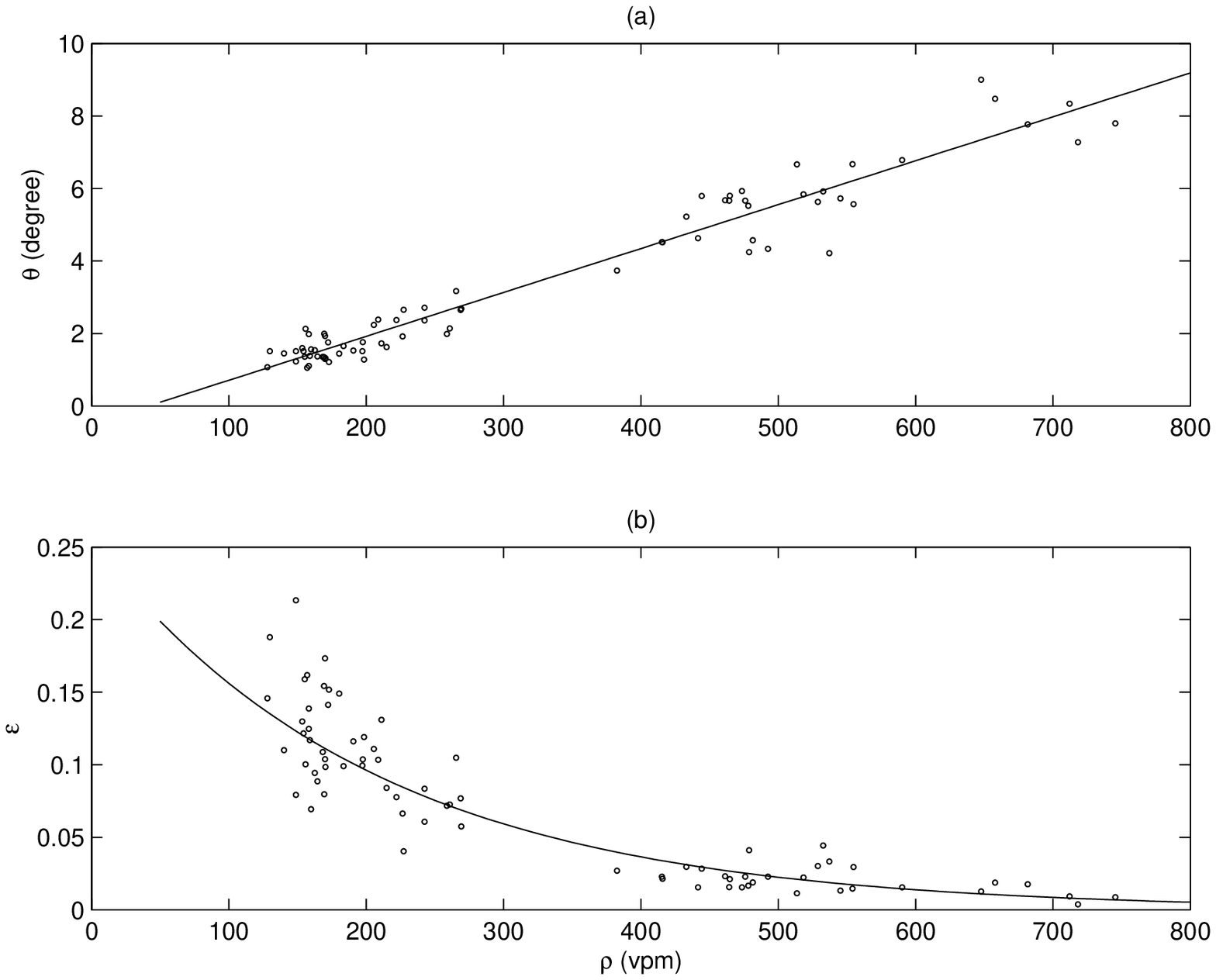}\caption{Characteristics of lane-changes at different densities with a threshold of $w$}\label{lc20081207characteristics}
\ec\efg
\bfg\bc
\includegraphics[height=4in]{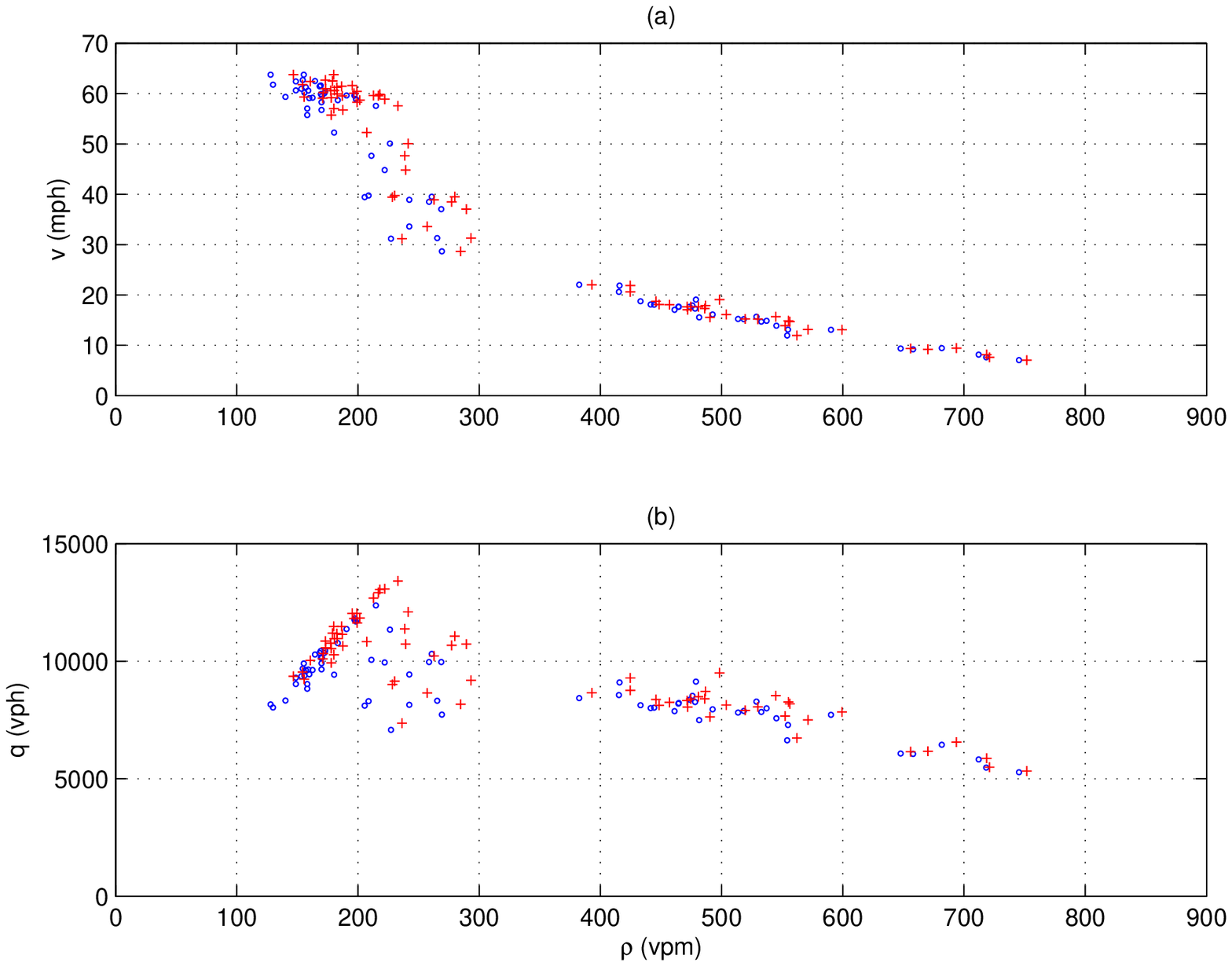}\caption{Fundamental diagrams with a threshold of $w$: Blue dots with lane-changing effect, and red dots without lane-changing effect}\label{lc20081207fd}
\ec\efg

With a lane-changing threshold of $\Delta y=1.5w$, from \reff{lc20081207characteristics2}, we find a linear relationship between density and lane-changing angle
\bqs
\theta&=&-0.3222+0.0086 \r,
\eqs
with R-square 0.9139.
We can also find the following relationship between density and $\epsilon$
\bqs
\epsilon&=&0.5579 e^{-0.0048 \r},
\eqs
with R-square 0.8917. 
From \reff{lc20081207fd2}, we can see that $q=\r V(\r)$ has a capacity of 15014 mph when $\r=251$ vpm, and $q=\r V((1+\epsilon) \r)$ has a capacity of 12369 mph when $\r=215$ vpm. The lane-changes cause $17.62\%$ capacity reduction. 
Potential capacity drop is more than twice than that for $\Delta y=w$. This is consistent with the observation in \reff{lc20081205trajectory}, since $t_{LC}$ is more than doubled with $\Delta y=1.5w$.

\bfg\bc
\includegraphics[height=4in]{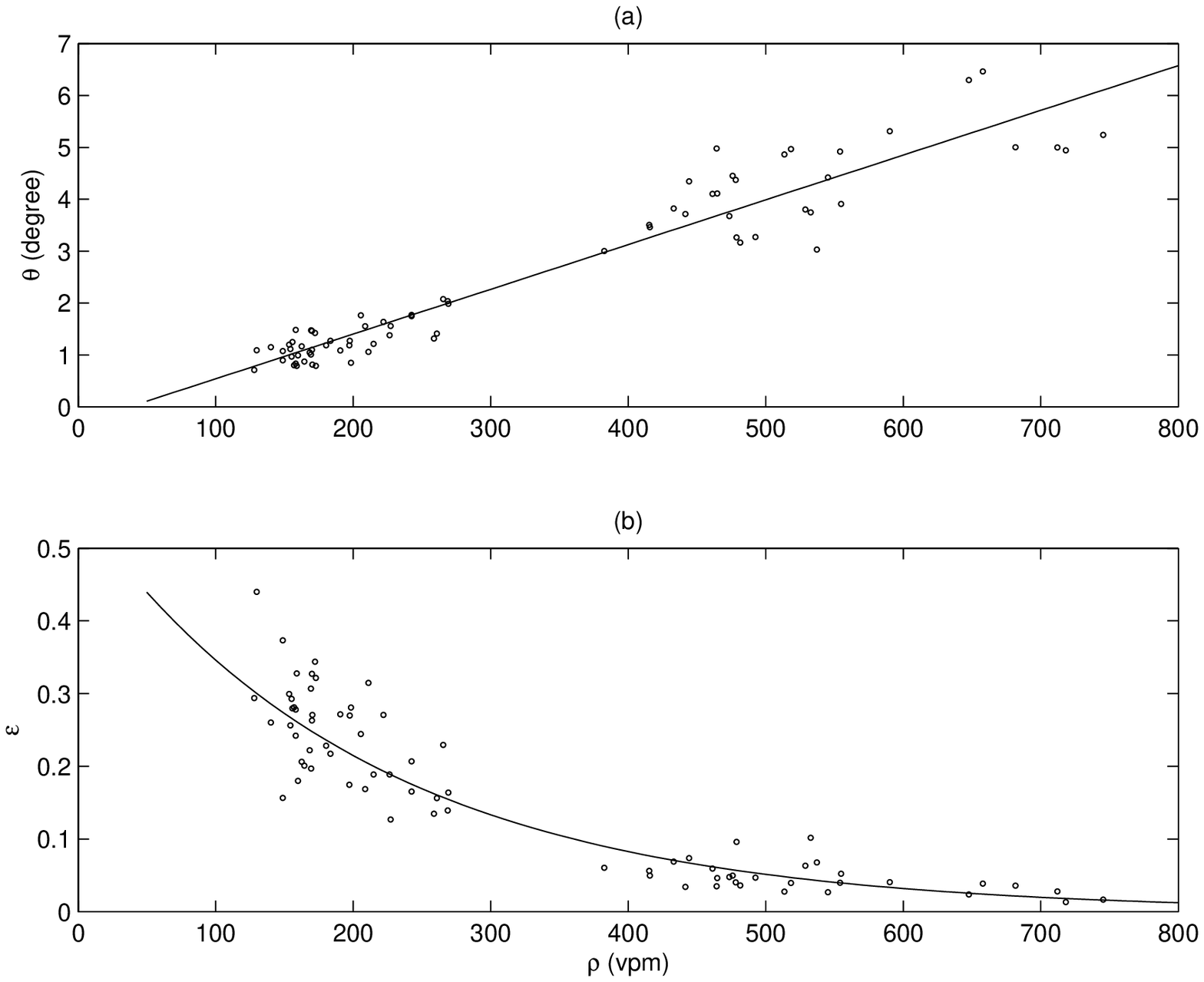}\caption{Characteristics of lane-changes at different densities with a threshold of $\frac 32 w$}\label{lc20081207characteristics2}
\ec\efg
\bfg\bc
\includegraphics[height=4in]{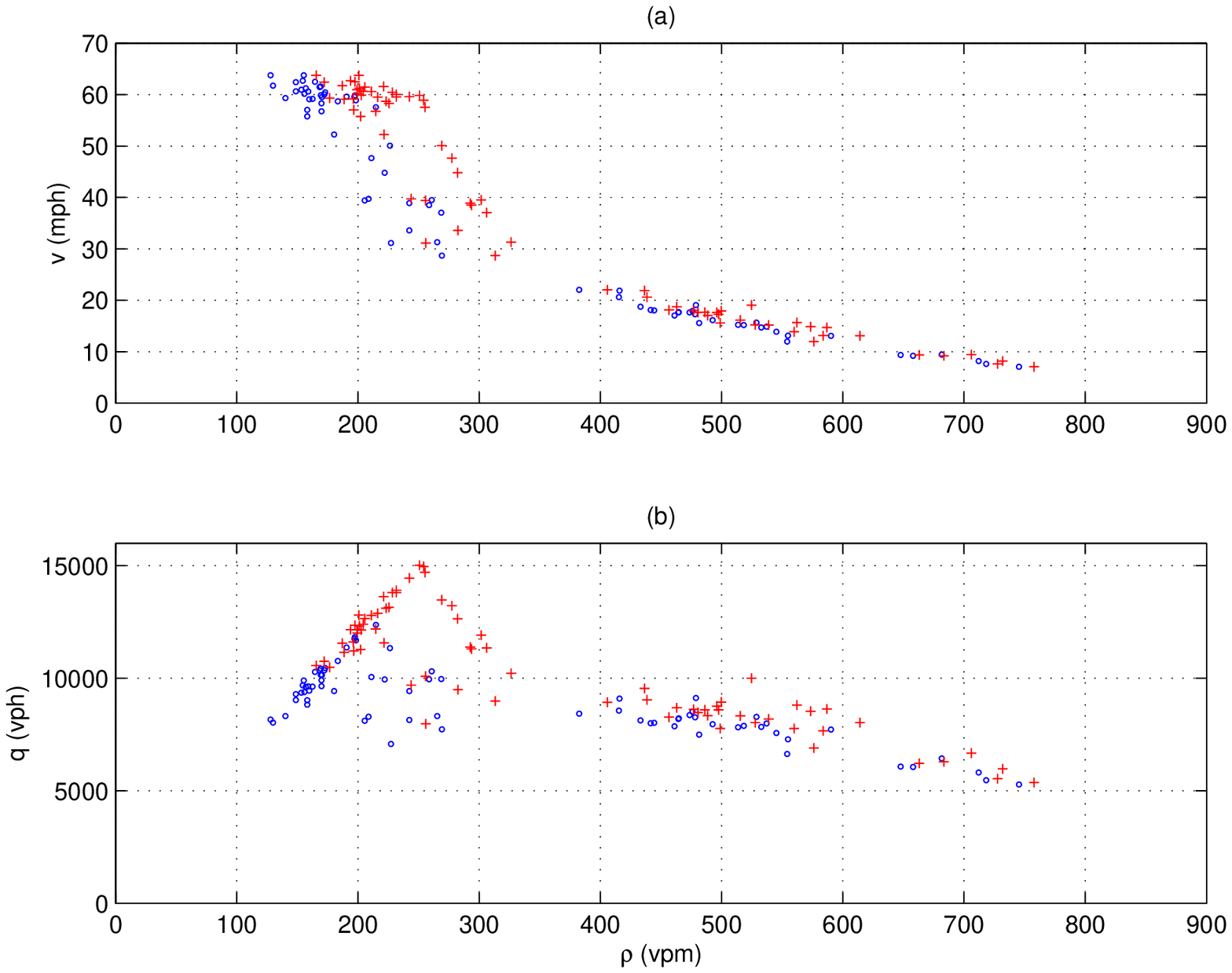}\caption{Fundamental diagrams with a threshold of $\frac 32 w$: Blue dots with lane-changing effect, and red dots without lane-changing effect}\label{lc20081207fd2}
\ec\efg

For the whole weaving location, we have density-dependent lane-changing intensity variables. More detailed analysis show that $\epsilon$ is also location-dependent \citep{jin2010lc}. Note that the kinematic wave model \refe{kw-weaving} is developed for a roadway with both lane-changing and normal sections. If we can obtain an $\epsilon-\r$ relation everywhere, \refe{kw-weaving} can be re-written as
\bqn
\rho_t+\left(\rho V(\rho(1+\epsilon(x,\rho)))\right)_x&=&0.
\eqn
If we define $Q(x,\rho)=\rho V(\rho(1+\epsilon(x,\rho)))$, this model is equivalent to an inhomogeneous LWR model. In this case, if we are able to introduce an inhomogeneity factor and write the equation above as a system of hyperbolic conservation laws, the analysis in Section 3 will apply for this general model. When we are not able to introduce an inhomogeneity factor, we can use the supply-demand method in \refe{s-d} for numerical solutions of the kinematic wave model of lane-changing traffic \refe{kw-weaving}.

\section{Conclusion}
In this paper, we presented a model for lane-changing traffic dynamics in the framework of kinematic waves. A key assumption is that the disruption lane-changing effect can be modeled by a modified speed-density relation with a new lane-changing intensity variable. 
Then, for location-dependent lane-changing intensities, we studies kinematic waves arising in lane-changing traffic and corresponding new definitions of local traffic supply and demand for computing the flow passing through a boundary.
With both theoretical derivations and observed data, we demonstrated that lane-changing intensities are highly related to road geometry, location, on-ramp/off-ramp traffic, lane-changing time,  vehicle speeds, and other traffic conditions. In particular, with 75 minutes of vehicle trajectories collected for a freeway section on interstate 80, we calibrated lane-changing intensities and corresponding fundamental diagrams for different lane-changing thresholds. The results support a functional relationship between density and the lane-changing effect for the whole weaving section, and the lane-changing angles are highly related to traffic density. Furthermore, it was suggested that lane-changes could cause 8\% to 18\% reduction in capacity, depending on the definition of a lane-changing threshold.

By incorporating lane-changing effects in the fundamental diagram, this study provides a simple framework to look into
lane-changing traffic dynamics in the framework of kinematic wave theory.
From this study, we can see that lane-changing traffic could cause capacity drop, different observed jam densities,
and fundamental diagram of reverse-$\lambda$ shape.  
Lane-changing traffic can also affect the formation and dissipation of shock and rarefaction waves on a roadway.  This simple model can be easily integrated into a
commodity-based kinematic wave simulation model of network traffic in \citep{jin2004network,jin2005fifo} to study
system-wide traffic dynamics. 

The model studied here bears certain limitations by assuming full balance among lanes and constant $\epsilon$.  Thus our model
is not intended for ramp weaves (HCM 2000\nocite{hcm2000}), but for usually balanced areas near a lane-drop or merging junction. 
By assuming an equilibrium speed-density relationship, \refe{kw-weaving} is not able to described detailed traffic dynamics when vehicles are accelerating or decelerating. 
In this study, we only studied the bottleneck effects of lane-changing traffic. In relatively sparse traffic, however, lane-changes could actually benefit the overall traffic flow due to their balancing effects, and such effects have yet to be included in $\epsilon$.
In addition, the new model considers the lateral interactions between vehicles in lane-changing areas at the aggregate level, and more detailed analysis of microscopic lane-changing decisions and maneuvers could give us more insights on the aggregate effects of lane-changing traffic.

With observed data, we calibrated a simple relationship between lane-changing intensity and traffic density. As shown in \refe{approxe}, lane-changing intensity is also related to on-ramp and off-ramp flows. For the weaving section in \reff{i80studyarea}, however, we do not have off-ramp flows in data sets 2 to 4 and therefore cannot study such relationships. In addition, it was shown that $\Delta y$ could significantly impact lane-changing intensities and therefore capacity reductions, and it can be calibrated by comparing road capacities with lane-changes and those without lane-changes.  For example, for the weaving section studied in Section 4, capacity reductions are about 8\% and 18\% when $\Delta y$  is set to be a vehicle¡¯s width and 1.5 times of a vehicle¡¯s width, respectively. In the future, we will be interested in calibrating the lane-changing threshold with more data sets or more detailed analysis.
In the future, we will also be interested in studying lane-changing intensities for different lane-changing areas, on-ramp/off-ramp flows, and other factors. Such a study would be helpful for understanding traffic macroscopic traffic dynamics in lane-changing areas and for developing possible ramp-metering and lane management strategies for improving the overall traffic flow.

\section*{Acknowledgments}
We would like to thank Dr. Shin-Ting (Cindy) Jeng for her valuable comments and discussions. The comments by several anonymous reviewers have been very helpful for improving the presentation of the paper. The views and results herein are the author's alone.

\section*{Appendix: Kinematic wave solutions of Types 2 to 10}
\bi
\item[Type 2] When $U_R$ is in Region 2, $Q(U_R)> Q(U_L)$. The Riemann problem is solved by a combination of a standing wave and a
forward traveling shock wave, with an intermediate state $U_1$ as shown in \reff{lctype2}. In this case,
$q=Q(U_L)$.
\bfg
\bc \includegraphics[height=8cm]{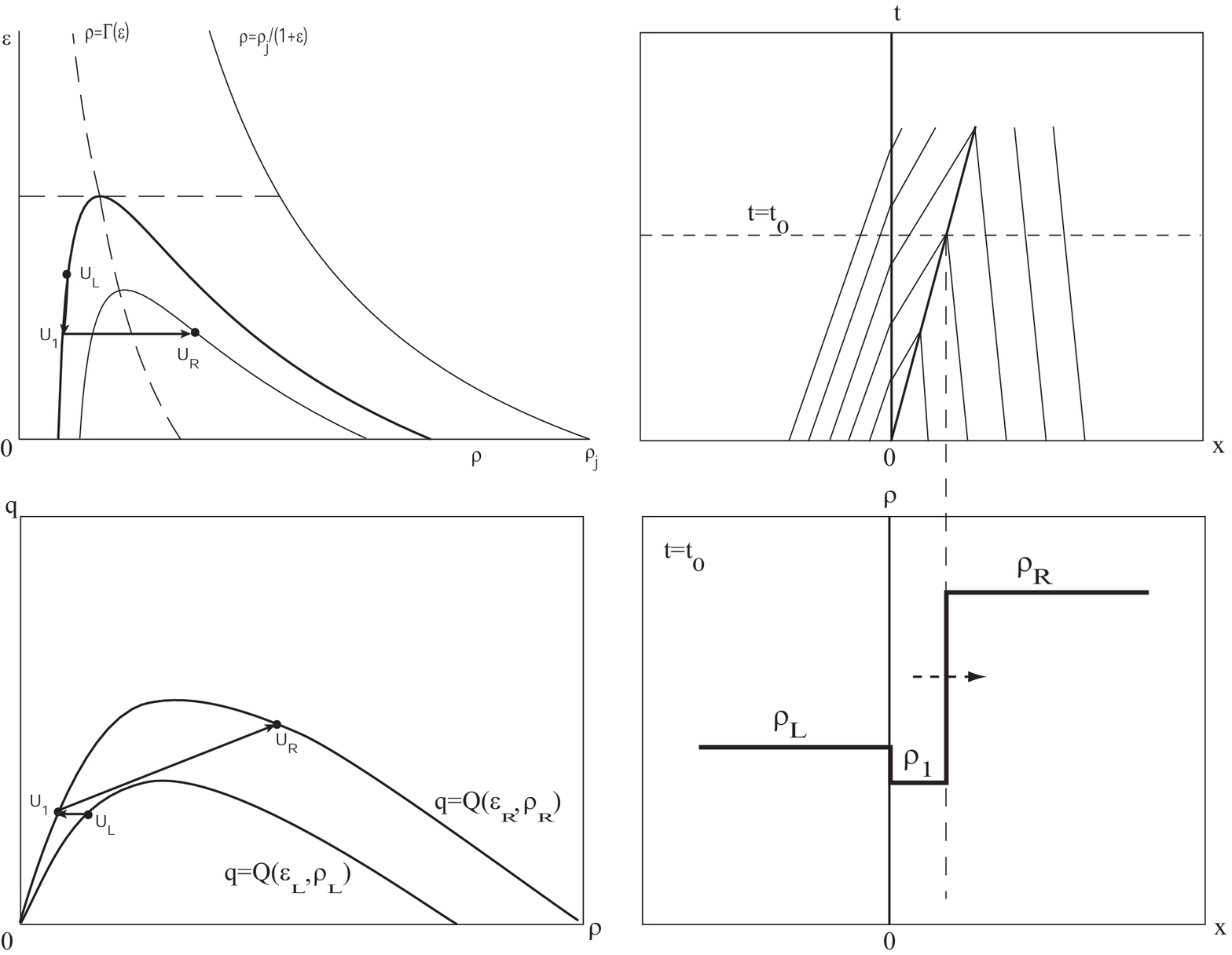}
\caption{An example for wave solutions to Riemann problem of Type 2}\label{lctype2}\ec
\efg

\item[Type 3] When $U_R$ is in Region 3, $\r_R> \Gamma(\e_R)$ and $Q(U_R)\leq Q(U_L)$. That is, $U_R$ is OC, and $Q(U_R)$ is not greater than $Q(U_L)$. The Riemann problem is solved by a combination of a backward traveling
shock wave and a standing wave, with an intermediate state $U_1$ as shown in \reff{lctype3}. In this case,
$q=Q(U_R)$.
\bfg
\bc \includegraphics[height=8cm]{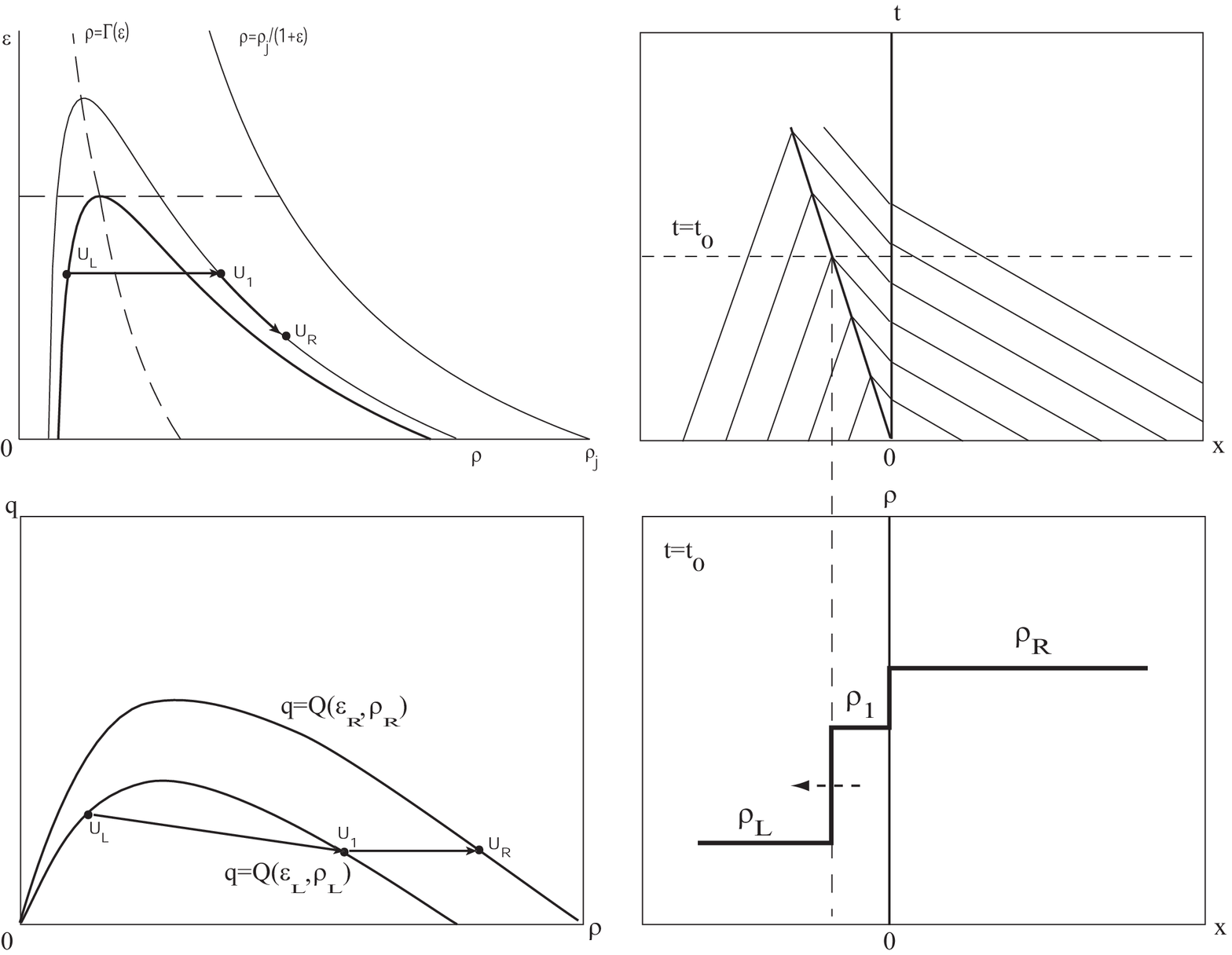}
\caption{An example for wave solutions to Riemann problem of Type 3}\label{lctype3}\ec
\efg
\item[Type 4] When $U_R$ is in Region 4, $\r_R\leq \Gamma(\e_R)$, $Q(U_R)\leq Q(U_L)$, and $Q(\e_R, \Gamma(\e_R))< Q(U_L)$. That is, $U_R$ is UC, and the capacity at $\e_R$ is smaller than $Q(U_L)$. The Riemann problem is solved by a combination of a backward traveling
shock wave, a standing wave, and a forward traveling rarefaction wave, with two intermediate states $U_1$ and
$U_2$ as shown in \reff{lctype4}, where $U_2=(\e_R,\Gamma(\e_R))$. In this case, $q=Q(\e_R,\Gamma(\e_R))$, which
is the capacity flow associated with $\e_R$.
\bfg
\bc \includegraphics[height=8cm]{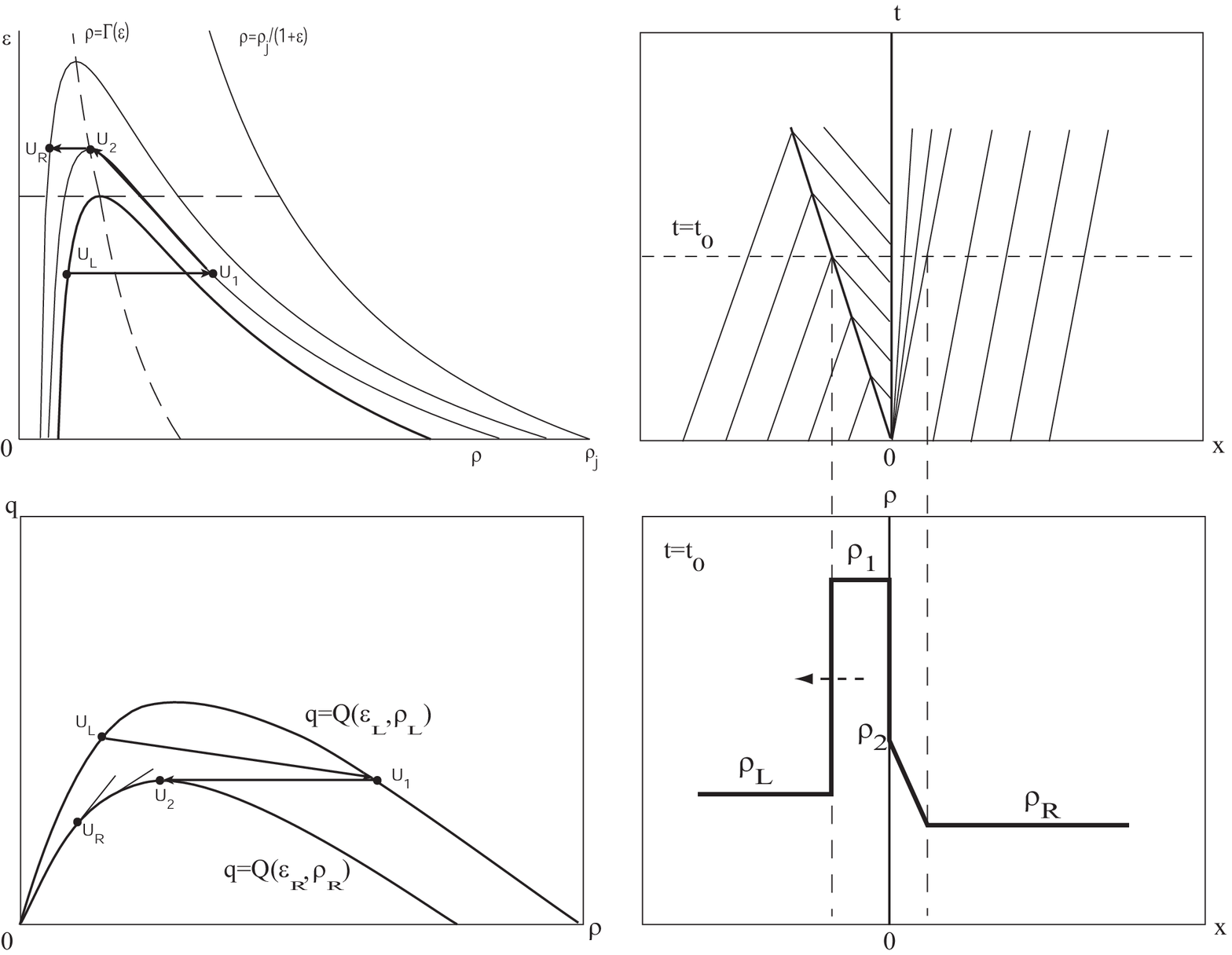}
\caption{An example for wave solutions to Riemann problem of Type 4}\label{lctype4}\ec
\efg

\item[Type 5] When $U_R$ is in Region 5, $\r_R\geq \Gamma(\e_R)$ and $Q(U_R)\leq Q(U_L)$. That is, $U_R$ is OC, and $Q(U_R)$ is not greater than $Q(U_L)$. The Riemann problem is solved by a combination of a backward traveling
shock wave and  a standing wave, with an intermediate state $U_1$ as shown in \reff{lctype5}. In this case,
$q=Q(U_R)$.
\bfg
\bc \includegraphics[height=8cm]{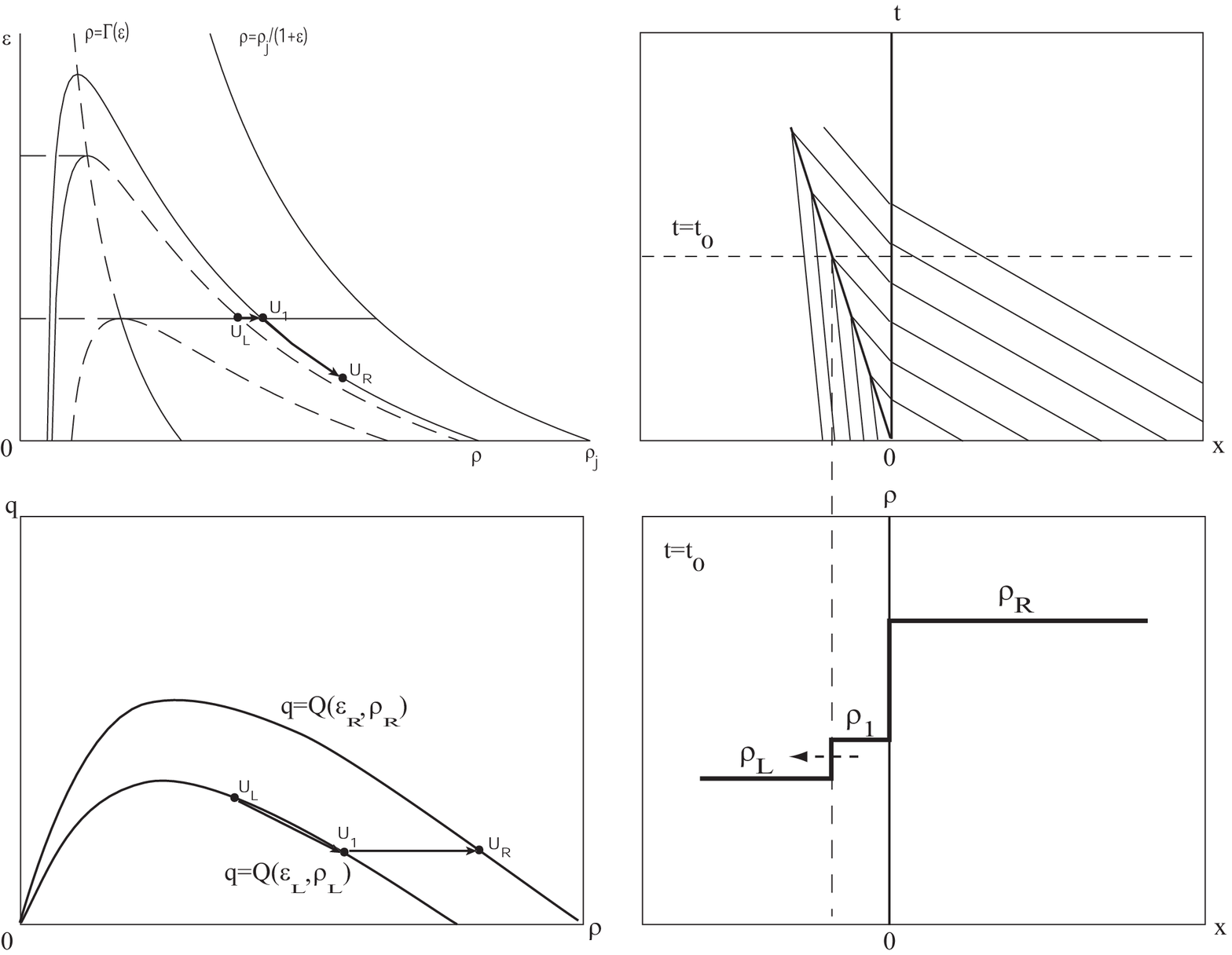}
\caption{An example for wave solutions to Riemann problem of Type 5}\label{lctype5}\ec
\efg
\item[Type 6] When $U_R$ is in Region 6, $\r_R\geq \Gamma(\e_R)$ and $Q(U_L)<Q(U_R)\leq Q(\e_L, \Gamma(\e_L))$. That is, $U_R$ is OC, and $Q(U_R)$ is between $Q(U_L)$ and the capacity at $\e_L$. The Riemann problem is solved by a combination of a backward traveling
rarefaction wave and a standing wave, with an intermediate state $U_1$ as shown in \reff{lctype6}. In this case,
$q=Q(U_R)$.
\bfg
\bc \includegraphics[height=8cm]{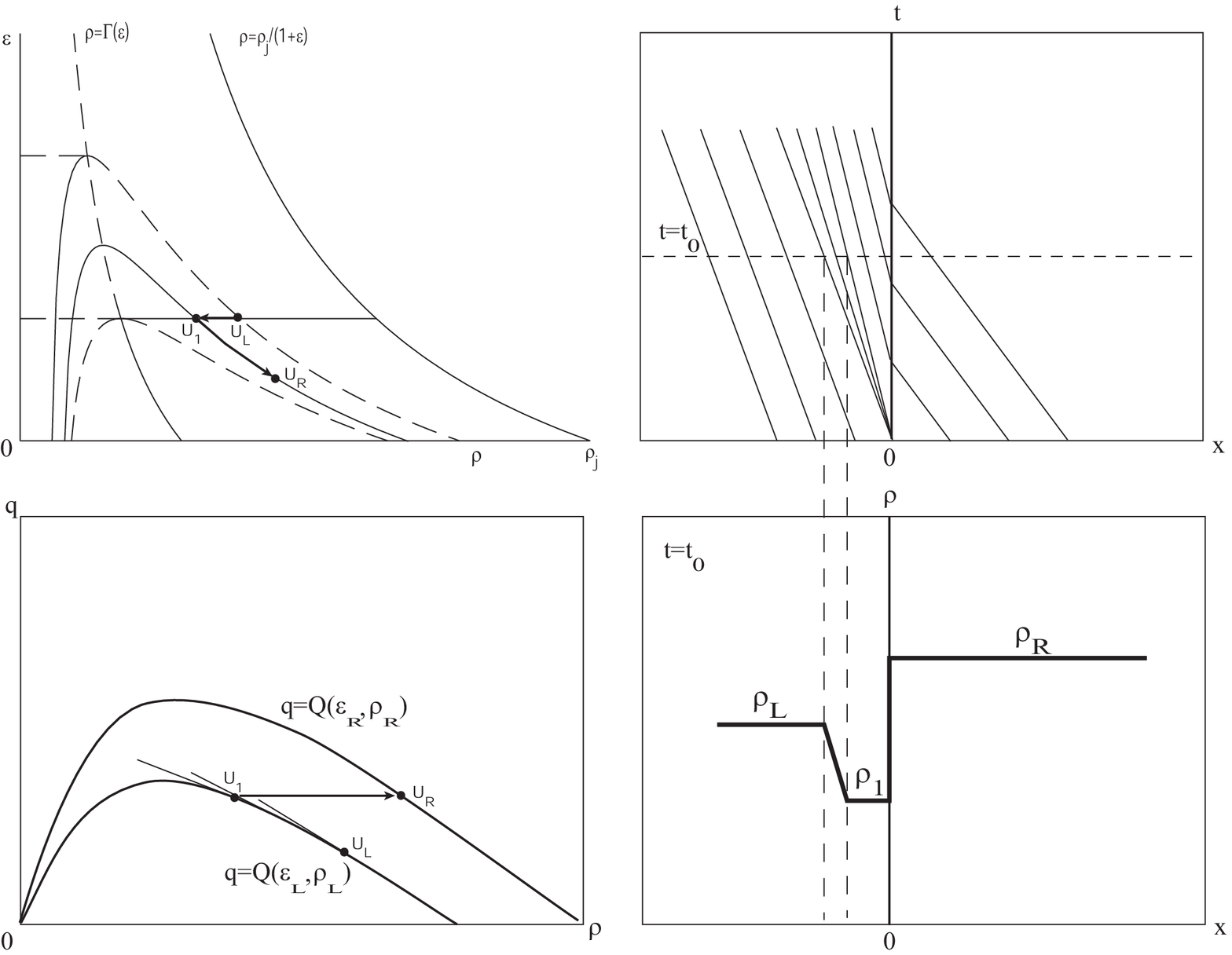}
\caption{An example for wave solutions to Riemann problem of Type 6}\label{lctype6}\ec
\efg
\item[Type 7] When $U_R$ is in Region 7, $\r_R< \Gamma(\e_R)$, $Q(U_R)\leq Q(\e_L, \Gamma(\e_L))$, but $Q(\e_R, \Gamma(\e_R))\geq Q(\e_L, \Gamma(\e_L))$. That is, $U_R$ is UC, and the capacity at $\e_L$ is between $Q(U_R)$ and the capacity at $\e_R$. The Riemann problem is solved by a combination of a backward traveling
rarefaction wave, a standing wave, and a forward traveling rarefaction wave, with two intermediate state $U_1$
and $U_2$ as shown in \reff{lctype7}, where $U_1=(\e_L, \Gamma(\e_L))$. In this case, $q=Q(U_2)=Q(U_1)=Q(\e_L,
\Gamma(\e_L))$, the capacity flow when $\e=\e_L$.
\bfg
\bc \includegraphics[height=8cm]{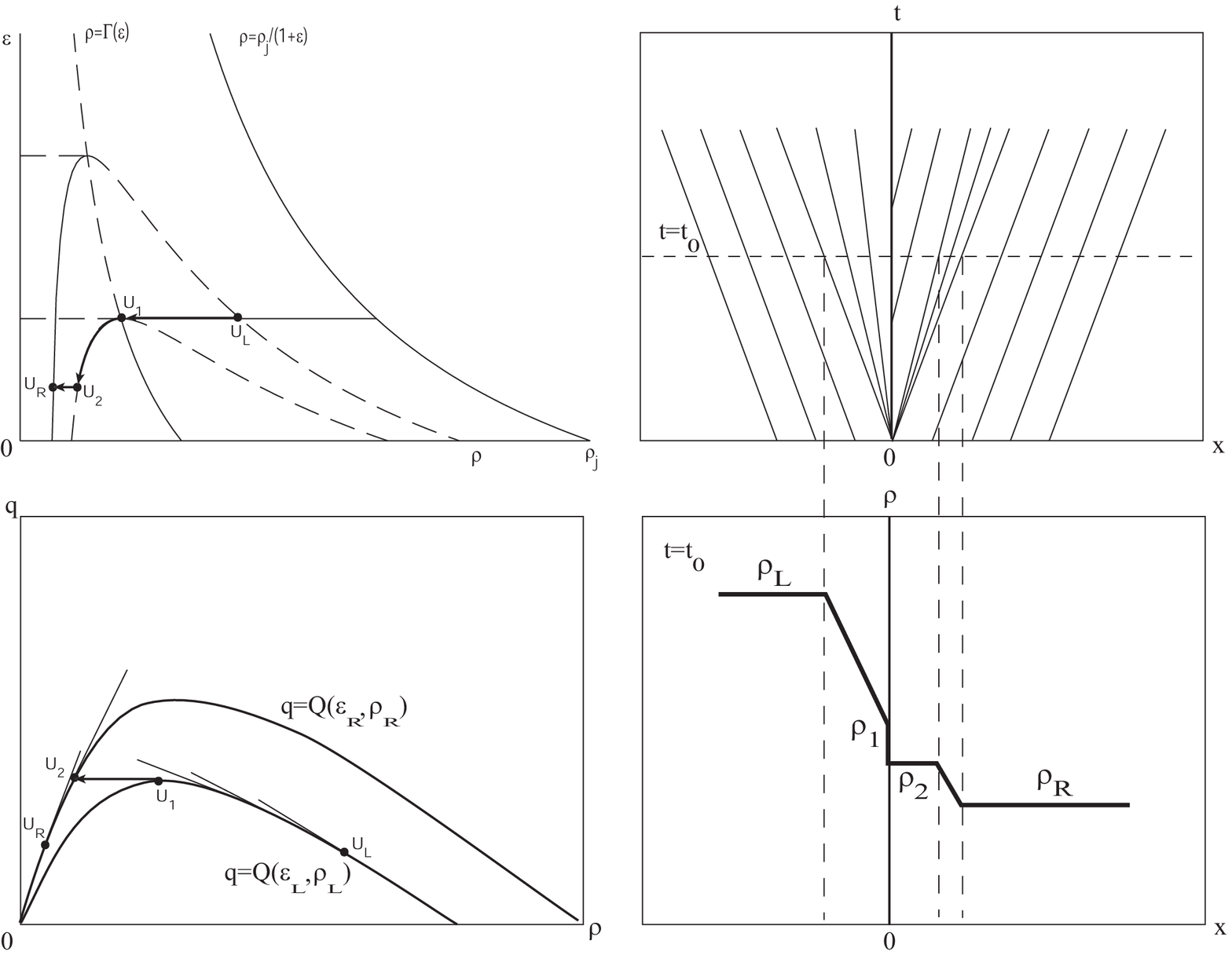}
\caption{An example for wave solutions to Riemann problem of Type 7}\label{lctype7}\ec
\efg
\item[Type 8] When $U_R$ is in Region 8, $Q(U_R)>Q(\e_L, \Gamma(\e_L))$. The Riemann problem is solved by a combination of a backward traveling
rarefaction wave, a standing wave, and a forward traveling shock wave, with two intermediate state $U_1$ and
$U_2$ as shown in \reff{lctype8}, where $U_1=(\e_L, \Gamma(\e_L))$. In this case, $q=Q(U_2)=Q(U_1)=Q(\e_L,
\Gamma(\e_L))$.
\bfg
\bc \includegraphics[height=8cm]{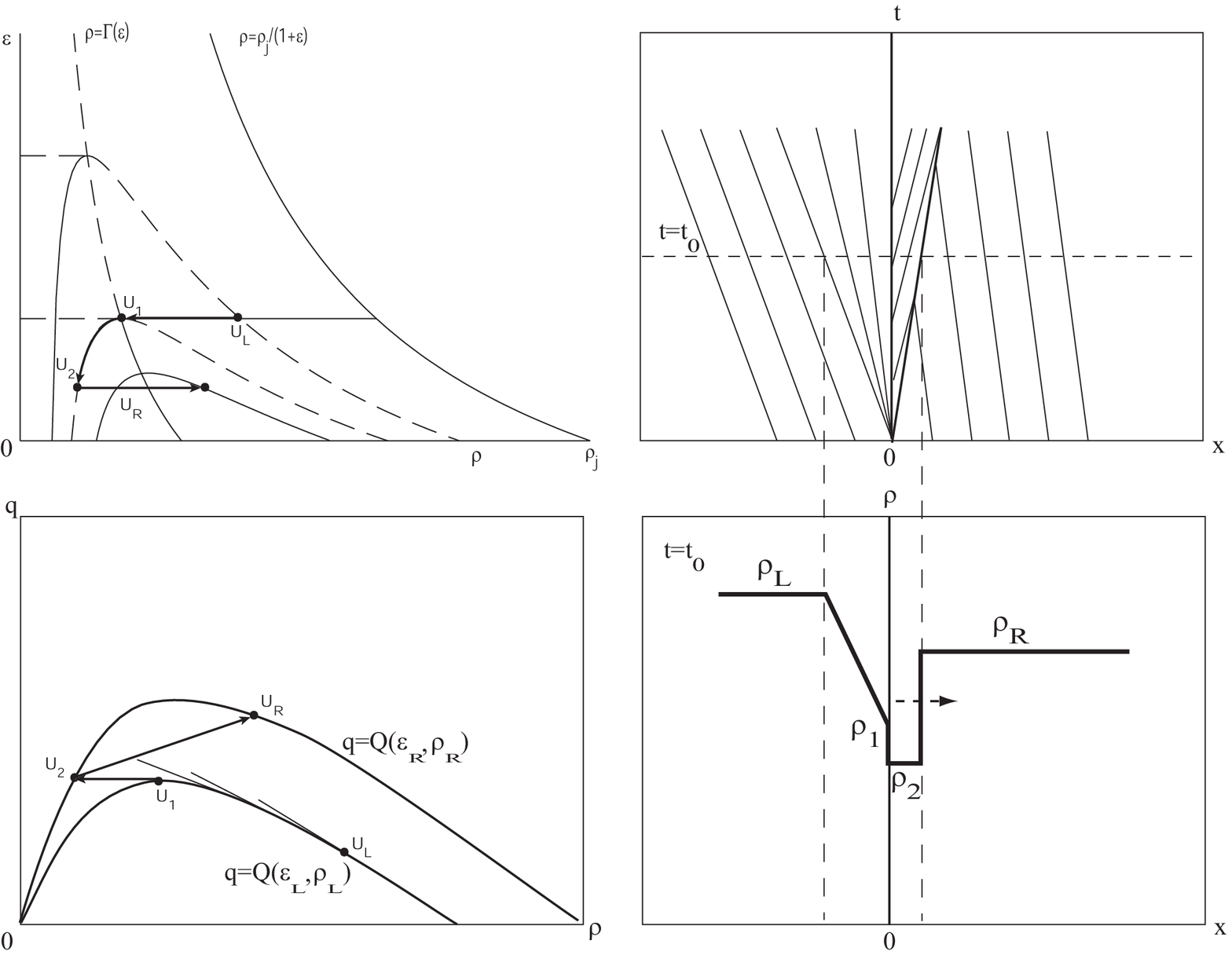}
\caption{An example for wave solutions to Riemann problem of Type 8}\label{lctype8}\ec
\efg
\item[Type 9] When $U_R$ is in Region 9, $\r_R<\Gamma(\e_R)$, $Q(U_R)<Q(U_L)$, and $Q(\e_R, \Gamma(\e_R))<Q(U_L)$. That is, $U_R$ is UC, and both $Q(U_R)$ and the capacity at $\e_R$ are smaller than $Q(U_L)$. The Riemann problem is solved by a combination of a backward traveling
shock wave, a standing wave, and a forward traveling rarefaction wave, with two intermediate state $U_1$ and
$U_2$ as shown in \reff{lctype9}, where $U_2=(\e_R, \Gamma(\e_R))$. In this case,
$q=Q(U_1)=Q(U_2)=Q(\e_R,\Gamma(\e_R))$.
\bfg
\bc \includegraphics[height=8cm]{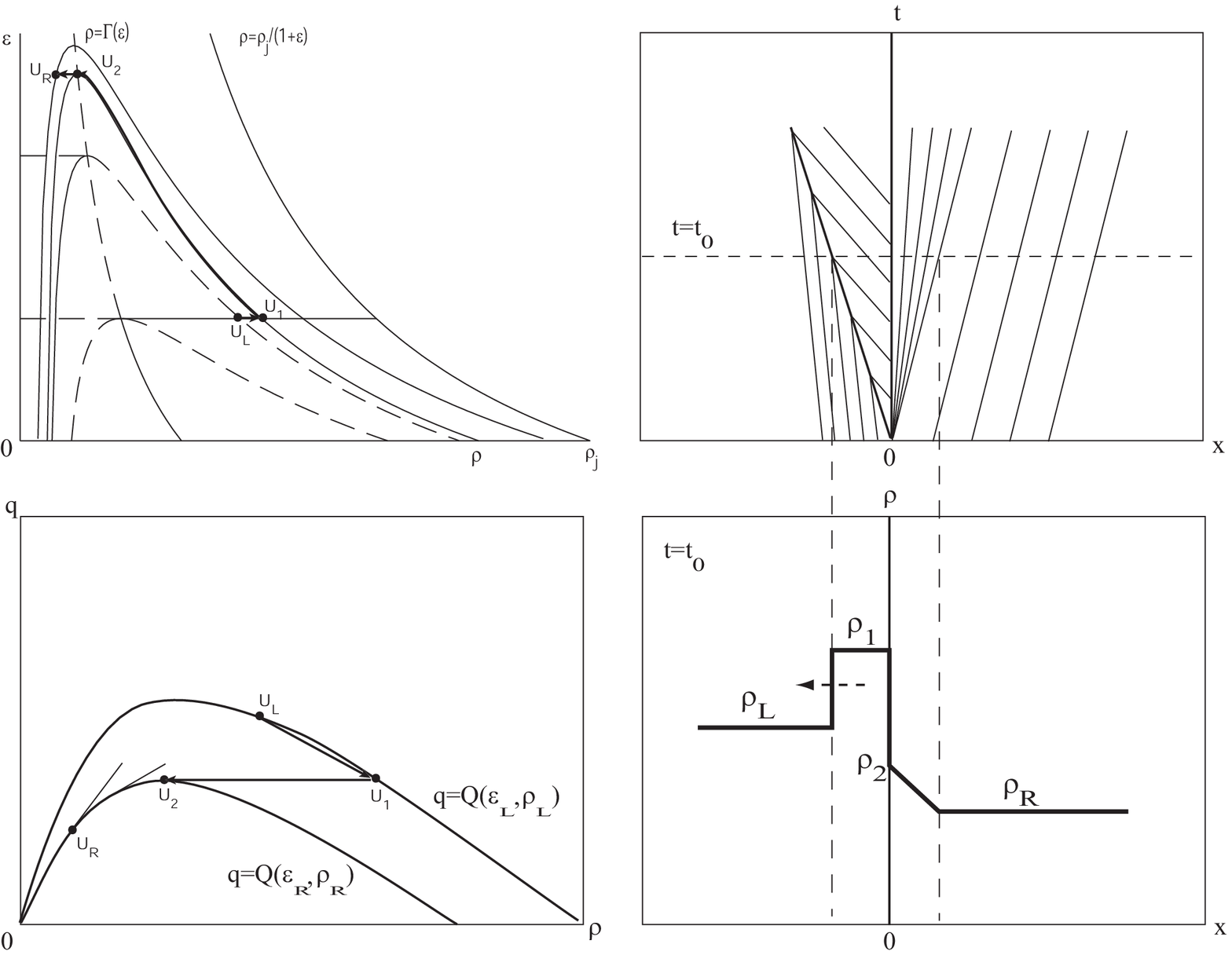}
\caption{An example for wave solutions to Riemann problem of Type 9}\label{lctype9}\ec
\efg
\item[Type 10] When $U_R$ is in Region 10, $\r_R<\Gamma(\e_R)$, and $Q(U_L)\leq Q(\e_R, \Gamma(\e_R))< Q(\e_L,\Gamma(\r_L))$. That is, $U_R$ is UC, and the capacity at $\e_R$ is between $Q(U_L)$ and the capacity at $\e_L$. The Riemann problem is solved by a combination of a backward traveling
rarefaction wave, a standing wave, and a forward traveling rarefaction wave, with two intermediate state $U_1$
and $U_2$ as shown in \reff{lctypea}, where $U_2=(\e_R, \Gamma(\e_R))$. In this case,
$q=Q(U_1)=Q(U_2)=Q(\e_R,\Gamma(\e_R))$.
\bfg
\bc \includegraphics[height=8cm]{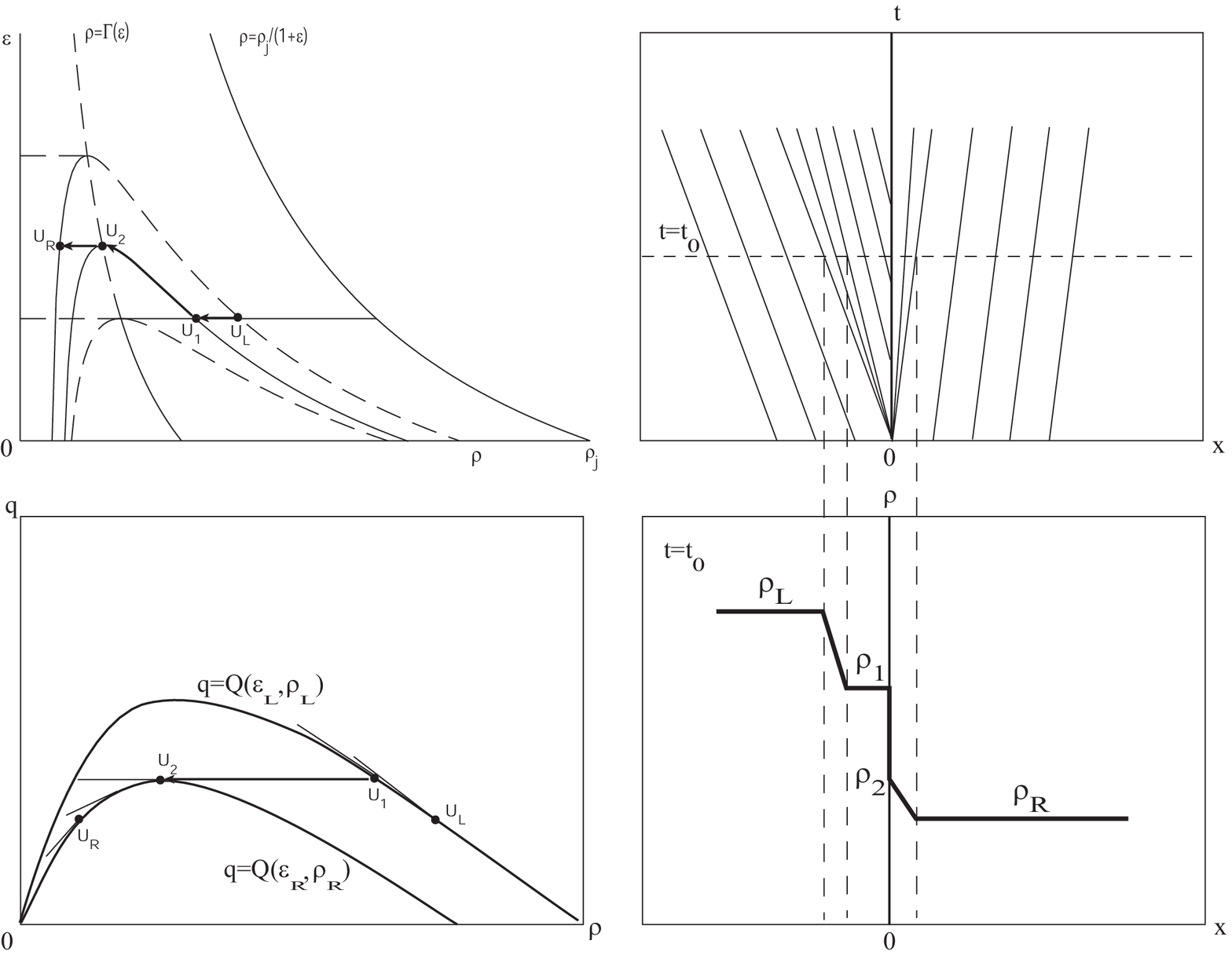}
\caption{An example for wave solutions to Riemann problem of Type 10}\label{lctypea}\ec
\efg
\ei

\begin{thebibliography}{61}
\expandafter\ifx\csname natexlab\endcsname\relax\def\natexlab#1{#1}\fi
\expandafter\ifx\csname url\endcsname\relax
  \def\url#1{\texttt{#1}}\fi
\expandafter\ifx\csname urlprefix\endcsname\relax\fi

\bibitem[{{Cambridge Systematics, Inc.}(2004)}]{ngsim2004data235}
{Cambridge Systematics, Inc.}, September 2004. {NGSIM BHL Data Analysis}. Tech.
  rep., summary Report, Prepared for Federal Highway Administration.

\bibitem[{{Cambridge Systematics, Inc.}(2005{\natexlab{a}})}]{ngsim2005data400}
{Cambridge Systematics, Inc.}, September 2005{\natexlab{a}}. {NGSIM I-80 Data
  Analysis (4:00 p.m. to 4:15 p.m.)}. Tech. rep., summary Report, Prepared for
  Federal Highway Administration.

\bibitem[{{Cambridge Systematics, Inc.}(2005{\natexlab{b}})}]{ngsim2005data500}
{Cambridge Systematics, Inc.}, September 2005{\natexlab{b}}. {NGSIM I-80 Data
  Analysis (5:00 p.m. to 5:15 p.m.)}. Tech. rep., summary Report, Prepared for
  Federal Highway Administration.

\bibitem[{{Cambridge Systematics, Inc.}(2005{\natexlab{c}})}]{ngsim2005data515}
{Cambridge Systematics, Inc.}, September 2005{\natexlab{c}}. {NGSIM I-80 Data
  Analysis (5:15 p.m. to 5:30 p.m.)}. Tech. rep., summary Report, Prepared for
  Federal Highway Administration.

\bibitem[{Cassidy and Rudjanakanoknad(2005)}]{cassidy2005merge}
Cassidy, M., Rudjanakanoknad, J., 2005. {Increasing the capacity of an isolated
  merge by metering its on-ramp}. Transportation Research Part B 39~(10),
  896--913.

\bibitem[{Cassidy et~al.(1989)Cassidy, Skabardonis, and
  May}]{cassidy1989weaving}
Cassidy, M., Skabardonis, A., May, A.~D., 1989. Operation of major freeway
  weaving sections: recent empirical evidence. Transportation Research Record:
  Journal of the Transportation Research Board 1225, 61--72.

\bibitem[{Cassidy and Bertini(1999)}]{cassidy1999bottlenecks}
Cassidy, M.~J., Bertini, R.~L., 1999. Some traffic features at freeway
  bottlenecks. Transportation Research Part B 33~(1), 25--42.

\bibitem[{Cassidy and May(1991)}]{cassidy1991weaving}
Cassidy, M.~J., May, A.~D., 1991. Proposed analytical technique for estimating
  capacity and level of service of major freeway weaving sections.
  Transportation Research Record: Journal of the Transportation Research Board
  1320, 99--109.

\bibitem[{Chang and Kao(1991)}]{chang1991characteristics}
Chang, G.-L., Kao, Y.-M., November 1991. An empirical investigation of
  macroscopic lane-changing characteristics on uncongested multilane freeways.
  Transportation Research Part A 25~(6), 375--389.

\bibitem[{Coifman(2003)}]{Coifman2003inflow}
Coifman, B., 2003. Estimating density and lane inflow on a freeway segment.
  Transportation Research Part A 37~(8), 689--701.

\bibitem[{Daganzo(1995)}]{daganzo1995ctm}
Daganzo, C.~F., 1995. The cell transmission model \m{II}: Network traffic.
  Transportation Research Part B 29~(2), 79--93.

\bibitem[{Daganzo(2002)}]{daganzo2002behavior}
Daganzo, C.~F., 2002. A behavioral theory of multi-lane traffic flow. \m{Part
  I}: Long homogeneous freeway sections. \m{II}: Merges and the onset of
  congestion. Transportation Research Part B 36, 131--169.

\bibitem[{{Del Castillo} and Benitez(1995)}]{delcastillo1995fd_empirical}
{Del Castillo}, J.~M., Benitez, F.~G., 1995. On the functional form of the
  speed-density relationship - \m{II}: Empirical investigation. Transportation
  Research Part B 29~(5), 391--406.

\bibitem[{Eads et~al.(2000)Eads, Rouphail, May, and Hall}]{eads2000hcm}
Eads, B.~S., Rouphail, N.~M., May, A.~M., Hall, F., 2000. Freeway facility
  methodology in ``{H}ighway {C}apacity {M}anual" 2000. Transportation Research
  Record: Journal of the Transportation Research Board 1710, 171--180.

\bibitem[{Fazio and Rouphail(1986)}]{fazio1986weaving}
Fazio, J., Rouphail, N.~M., 1986. Freeway weaving sections: comparison and
  refinement of design and operations analysis procedures. Transportation
  Research Record: Journal of the Transportation Research Board 1091, 101--109.

\bibitem[{{Federal highway administration}(2000)}]{hcm2000}
{Federal highway administration}, 2000. {Highway capacity manual}. Tech. rep.,
  Transportation Research Board, National Research Council, Washington, D.C.

\bibitem[{{Federal Highway Administration}(2006)}]{fhwa2006ngsim}
{Federal Highway Administration}, December 2006. {Next Generation SIMulation
  Fact Sheet}. Tech. rep., {FHWA-HRT-06-135}.

\bibitem[{Fitzpatrick and Nowlin(1996)}]{fitzpatrick1996weaving}
Fitzpatrick, K., Nowlin, L., 1996. One-sided weaving operations on one-way
  frontage roads. Transportation Research Record: Journal of the Transportation
  Research Board 1555, 42--49.

\bibitem[{Gazis et~al.(1961)Gazis, Herman, and Rothery}]{gazis1961follow}
Gazis, D.~C., Herman, R., Rothery, R.~W., 1961. Nonlinear follow-the-leader
  models of traffic flow. Operations Research 9~(4), 545--567.

\bibitem[{Gazis et~al.(1962)Gazis, Herman, and Weiss}]{gazis1962multilane}
Gazis, D.~C., Herman, R., Weiss, G.~H., 1962. Density oscillations between
  lanes of a multilane highway. Operations Research 10~(5), 658--667.

\bibitem[{Gipps(1986)}]{gipps1986changing}
Gipps, P.~G., 1986. A model for the structure of lane changing decisions.
  Transportation Research Part B 20~(5), 403--414.

\bibitem[{Godunov(1959)}]{godunov1959}
Godunov, S.~K., 1959. A difference method for numerical calculations of
  discontinuous solutions of the equations of hydrodynamics. Matematicheskii
  Sbornik 47, 271--306, in Russian.

\bibitem[{Golob et~al.(2004)Golob, Recker, and Alvarez}]{golob2004safety}
Golob, T.~F., Recker, W.~W., Alvarez, V.~M., January 2004. Safety aspects of
  freeway weaving sections. Transportation Research Part A: Policy and Practice
  38~(1), 35--51.

\bibitem[{Greenshields(1935)}]{greenshields1935capacity}
Greenshields, B.~D., 1935. A study in highway capacity. Highway Research Board
  Proceedings 14, 448--477.

\bibitem[{Haberman(1977)}]{haberman1977model}
Haberman, R., 1977. Mathematical models. Prentice Hall, Englewood Cliffs, NJ.

\bibitem[{Hall and Agyemang-Duah(1991)}]{hall1991capacity}
Hall, F.~L., Agyemang-Duah, K., 1991. Freeway capacity drop and the definition
  of capacity. Transportation Research Record: Journal of the Transportation
  Research Board 1320, 91--98.

\bibitem[{Holland and Woods(1997)}]{holland1997continuum}
Holland, E.~N., Woods, A.~W., November 1997. A continuum model for the
  dispersion of traffic on two-lane roads. Transportation Research Part B
  31~(6), 473--485.

\bibitem[{Isaacson and Temple(1992)}]{isaacson1992resonance}
Isaacson, E.~I., Temple, J.~B., 1992. Nonlinear resonance in systems of
  conservation laws. SIAM Journal on Applied Mathematics 52~(5), 1260--1278.

\bibitem[{Jin(2010)}]{jin2010lc}
Jin, W.-L., 2010. {Macroscopic characteristics of lane-changing vehicular
  traffic}. In: Proceedings of Transportation Research Board Annual Meeting. To
  be presented.

\bibitem[{Jin et~al.(2009)Jin, Chen, and Puckett}]{jin2009sd}
Jin, W.-L., Chen, L., Puckett, E.~G., 2009. {Supply-demand diagrams and a new
  framework for analyzing the inhomogeneous Lighthill-Whitham-Richards model}.
  {Proceedings of the 18th International Symposium on Transportation and
  Traffic Theory (ISTTT18)}, 603--635.

\bibitem[{Jin and Jayakrishnan(2005)}]{jin2005fifo}
Jin, W.-L., Jayakrishnan, R., 2005. First-in-first-out properties of a
  commodity-based kinematic wave simulation model. Transportation Research
  Record: Journal of the Transportation Research Board 1934, 197--207.

\bibitem[{Jin and Zhang(2003)}]{jin2003inhlwr}
Jin, W.-L., Zhang, H.~M., August 2003. The inhomogeneous kinematic wave traffic
  flow model as a resonant nonlinear system. Transportation Science 37~(3),
  294--311.

\bibitem[{Jin and Zhang(2004)}]{jin2004network}
Jin, W.-L., Zhang, H.~M., 2004. A multicommodity kinematic wave simulation
  model of network traffic flow. Transportation Research Record: Journal of the
  Transportation Research Board 1883, 59--67.

\bibitem[{Jin et~al.(2006)Jin, Zhang, and Chu}]{jin2006fifo}
Jin, W.-L., Zhang, Y., Chu, L., 2006. Measuring first-in-first-out violation
  among vehicles. In: Proceedings of Transportation Research Board Annual
  Meeting.

\bibitem[{Kesting et~al.(2007)Kesting, Treiber, and Helbing}]{kesting2007lc}
Kesting, A., Treiber, M., Helbing, D., 2007. {General Lane-Changing Model MOBIL
  for Car-Following Models}. Transportation Research Record: Journal of the
  Transportation Research Board 1999~(-1), 86--94.

\bibitem[{Klar and Wegener(1999)}]{Klar1999multilane}
Klar, A., Wegener, R., 1999. A hierarchy of models for multilane vehicular
  traffic {I}: Modeling. {II}: Numerical investigations. SIAM Journal on
  Applied Mathematics 59~(3), 983--1011.

\bibitem[{Koshi et~al.(1983)Koshi, Iwasaki, and Ohkura}]{koshi1983fd}
Koshi, M., Iwasaki, M., Ohkura, I., 1983. Some findings and an overview on
  vehicular flow characteristics. In: Hurdle, V.~F., Hauer, R., Stewart, G.~N.
  (Eds.), Proceedings of the Eighth International Symposium on Transportation
  and Traffic Theory. University of Toronto Press, Toronto, Ontario, pp.
  403--426.

\bibitem[{Laval and Daganzo(2006)}]{laval2006lc}
Laval, J., Daganzo, C., 2006. {Lane-changing in traffic streams}.
  Transportation Research Part B 40~(3), 251--264.

\bibitem[{Lax(1972)}]{lax1972shock}
Lax, P.~D., 1972. Hyperbolic systems of conservation laws and the mathematical
  theory of shock waves. SIAM, Philadelphia, Pennsylvania.

\bibitem[{Lebacque(1996)}]{lebacque1996godunov}
Lebacque, J.~P., 1996. The {G}odunov scheme and what it means for first order
  traffic flow models. In: The International Symposium on Transportation and
  Traffic Theory. Lyon, France.

\bibitem[{Leisch(1979)}]{leisch1979weaving}
Leisch, J.~E., 1979. A new technique for design and analysis of weaving
  sections on freeways. ITE Journal 49~(3), 26--29.

\bibitem[{Lighthill and Whitham(1955)}]{lighthill1955lwr}
Lighthill, M.~J., Whitham, G.~B., 1955. On kinematic waves: \m{II. A} theory of
  traffic flow on long crowded roads. Proceedings of the Royal Society of
  London A 229~(1178), 317--345.

\bibitem[{Michalopoulos et~al.(1984)Michalopoulos, Beskos, and
  Yamauchi}]{michalopoulos1984multilane}
Michalopoulos, P.~G., Beskos, D.~E., Yamauchi, Y., 1984. Multilane traffic flow
  dynamics: Some macroscopic considerations. Transportation Research Part B
  18~(4/5), 377--395.

\bibitem[{Milam and Choa(1998)}]{milam1998cloverleaf}
Milam, R.~T., Choa, F., 1998. A comparison of partial and full partial
  cloverleaf interchange operations using the {CORSIM} micro-simulation model.
  In: ITE District 6 Annual Conference. San Jose, California.

\bibitem[{Munjal et~al.(1971)Munjal, Hsu, and Lawrence}]{munjal1971multilane}
Munjal, P.~K., Hsu, Y.~S., Lawrence, R.~L., 1971. Analysis and validation of
  lane-drop effects of multilane freeways. Transportation Research 5, 257--266.

\bibitem[{Munjal and Pipes(1971)}]{munjal1971multilane1}
Munjal, P.~K., Pipes, L.~A., 1971. Propagation of on-ramp density waves on
  uniform unidirectional multilane freeways. Transportation Research 5,
  241--255.

\bibitem[{Nelson and Kumar(2004)}]{nelson2004model}
Nelson, P., Kumar, N., 2004. Point constriction, interface and boundary
  conditions for the kinematic-wave model. In: Proceedings of Transportation
  Research Board Annual Meeting.

\bibitem[{Newell(1961)}]{newell1961nonlinear}
Newell, G.~F., 1961. Nonlinear effects in the dynamics of car following.
  Operations Research 9~(2), 209.

\bibitem[{Newell(1993)}]{newell1993sim}
Newell, G.~F., 1993. A simplified theory of kinematic waves in highway traffic
  \m{I}: General theory. \m{II}: Queuing at freeway bottlenecks. \m{III}:
  Multi-destination flows. Transportation Research Part B 27~(4), 281--313.

\bibitem[{Ostrom et~al.(1993)Ostrom, Leiman, and May}]{ostrom1993weaving}
Ostrom, B., Leiman, L., May, A.~D., 1993. Suggested procedures for analyzing
  freeway weaving sections. Transportation Research Record: Journal of the
  Transportation Research Board 1398, 42--48.

\bibitem[{Pahl(1972)}]{Pahl1972lc}
Pahl, J., 1972. Lane-change frequencies in freeway traffic flow. Highway
  Research Record 409, 17--25.

\bibitem[{Richards(1956)}]{richards1956lwr}
Richards, P.~I., 1956. Shock waves on the highway. Operations Research 4,
  42--51.

\bibitem[{Roess et~al.(1974)Roess, McShane, and Pignataro}]{roess1974weaving}
Roess, R.~P., McShane, W.~R., Pignataro, L.~J., 1974. Configuration, design,
  and analysis of weaving sections. Transportation Research Record: Journal of
  the Transportation Research Board 489, 1--12.

\bibitem[{Sheu and Ritchie(2001)}]{sheu2001changing}
Sheu, J.-B., Ritchie, S.~G., 2001. Stochastic modeling and real-time prediction
  of vehicular lane-changing behavior. Transportation Research Part B 35~(7),
  695--716.

\bibitem[{Toledo et~al.(2003)Toledo, Koutsopoulos, and
  Ben-Akiva}]{toledo2003changing}
Toledo, T., Koutsopoulos, H.~N., Ben-Akiva, M.~E., 2003. Modeling integrated
  lane-changing behavior. In: Proceedings of Transportation Research Board
  Annual Meeting.

\bibitem[{Wang and Prevedouros(1998)}]{wangyh1998comparison}
Wang, Y., Prevedouros, P.~D., 1998. Comparison of {INTEGRATION}, {TSIS/CORSIM}
  and {WATSim} in replicating volumes and speeds on three small networks.
  Transportation Research Record: Journal of the Transportation Research Board
  1644, 80--92.

\bibitem[{Whitham(1974)}]{whitham1974PW}
Whitham, G.~B., 1974. Linear and nonlinear waves. John Wiley and Sons, New
  York.

\bibitem[{Windover and May(1994)}]{windover1994weaving}
Windover, J.~R., May, A.~D., 1994. Revisions to level {D} methodology of
  analyzing freeway ramp weaving sections. Transportation Research Record:
  Journal of the Transportation Research Board 1457, 43--49.

\bibitem[{Worrall and Bullen(1970)}]{Worrall1970lc1}
Worrall, R.~D., Bullen, A.~G.~R., 1970. An empirical lane-changing model on
  multilane highways. Highway Research Record 303, 30--43.

\bibitem[{Worrall et~al.(1970)Worrall, Bullen, and Gur}]{Worrall1970lc}
Worrall, R.~D., Bullen, A.~G.~R., Gur, Y., 1970. An elementary stochastic model
  of lane-changing on a multilane highway. Highway Research Record 308, 1--12.

\bibitem[{Yang(1997)}]{yang1997dissertation}
Yang, Q., 1997. A simulation laboratory for evaluation of dynamic traffic
  management systems. Ph.D. thesis, Massachusetts Institute of Technology,
  Cambridge, Massachusetts.

\end{thebibliography}
\end {document}